\pdfoutput=1
\RequirePackage{ifpdf}
\ifpdf 
\documentclass[pdftex]{sigma}
\else
\documentclass{sigma}
\fi

\usepackage[all]{xy}
\usepackage{tikz}
\DeclareMathOperator{\co}{co}
\DeclareMathOperator{\Span}{Span}

\def\sw#1{{\sb{(#1)}}}

\newcommand{\ls}[1]{\ell(#1)^{\langle 1\rangle}}

\newcommand{\rs}[1]{\ell(#1)^{\langle 2\rangle}}

\newcommand{\lls}[2]{\ell_{#1}(#2)^{\langle 1\rangle}}

\newcommand{\rrs}[2]{\ell_{#1}(#2)^{\langle 2\rangle}}

\newcommand{\qcp}[1]{{C(\mathbb{R} P^{#1}_\mathcal{T})}}

\newcommand{\proj}[2]{\mathbb{#1}P_\mathcal{T}^#2}

\begin{document}

\newcommand{\arXivNumber}{1403.7818}

\allowdisplaybreaks

\renewcommand{\thefootnote}{$\star$}

\renewcommand{\PaperNumber}{088}

\FirstPageHeading

\ShortArticleName{Piecewise Principal Coactions of Co-Commutative Hopf Algebras}

\ArticleName{Piecewise Principal Coactions\\
of Co-Commutative Hopf Algebras\footnote{This paper is a~contribution to the Special Issue on Noncommutative Geometry
and Quantum Groups in honor of Marc A.~Rief\/fel.
The full collection is available at
\href{http://www.emis.de/journals/SIGMA/Rieffel.html}{http://www.emis.de/journals/SIGMA/Rieffel.html}}}

\Author{Bartosz ZIELI\'NSKI}

\AuthorNameForHeading{B.~Zieli\'nski}

\Address{Department of Computer Science, Faculty of Physics and Applied Informatics,\\
University of \L{}\'od\'z, Pomorska 149/153 90-236 \L{}\'od\'z, Poland}
\Email{\href{mailto:bzielinski@uni.lodz.pl}{bzielinski@uni.lodz.pl}}

\ArticleDates{Received March 31, 2014, in f\/inal form August 11, 2014; Published online August 18, 2014}

\Abstract{Principal comodule algebras can be thought of as objects representing principal bundles in non-commutative
geometry.
A~crucial component of a~principal comodule algebra is a~strong connection map.
For some applications it suf\/f\/ices to prove that such a~map exists, but for others, such as computing the associated
bundle projectors or Chern--Galois characters, an explicit formula for a~strong connection is necessary.
It has been known for some time how to construct a~strong connection map on a~multi-pullback comodule algebra from
strong connections on multi-pullback components, but the known explicit general formula is unwieldy.
In this paper we derive a~much easier to use strong connection formula, which is not, however, completely general, but
is applicable only in the case when a~Hopf algebra is co-commutative.
Because certain linear splittings of projections in multi-pullback comodule algebras play a~crucial role in our
construction, we also devote a~signif\/icant part of the paper to the problem of existence and explicit formulas for such
splittings. Finally, we show example application of our work.}


\Classification{58B32; 16T05}

\renewcommand{\thefootnote}{\arabic{footnote}}
\setcounter{footnote}{0}

\vspace{-2mm}

\section{Introduction}

\looseness=-1
Let~$H$ be a~Hopf algebra (with bijective antipode), interpreted as a~Peter--Weyl algebra of functions on a~quantum
group.
Principal $H$-comodule algebras can be loosely viewed as the algebras of appropriate classes of functions on
(non-commutative) principal bundles (\cite{bhms06} makes the relationship explicit in the classical case).
A~crucial ingredient in the def\/inition of principal comodule algebra is a~so called strong connection map.
For some applications it suf\/f\/ices to prove that a~strong connection map exists, for instance when proving principality
of a~comodule algebra (see, e.g.,~\cite{reduct}).
Other applications (see, e.g.,~\cite{bf,hms03, hms2, haj-wagn}), such as computing the associated bundle projector or
Chern--Galois character~\cite{bh04}, call for an explicit formula for this map.

Piecewise principal comodule algebras~\cite{cm2002,hkmz11} is an interesting class of principal comodule algebras for
which a~fair amount of examples recently appeared in the literature (see,
e.g.,~\cite{bhms,bf14,cp14,index-pair,hms2,hms,Rennie,reduct,haj-wagn}).
They can be understood as being glued (constructed as a~multi-pullback) from simpler parts which are principal.
In~\cite{hkmz11} (cf.~the generalization in~\cite{haj-wagn}) it was proven that piecewise principal comodule algebras
are, in fact, principal.
The paper contains a~derivation of the explicit formula for a~strong connection on a~pullback of two principal
extensions from the ``local'' strong connections on pullback components and an appropriate choice of splittings of the
gluing maps.
If the piecewise comodule algebra is a~multipullback one can present this multipullback as an iterated pullback, and
then iterate the formula.
Unfortunately, in practice, already the second iteration of the formula from~\cite{hkmz11} becomes overly complicated.

\looseness=-1
In the paper we derive, under the assumption of the co-commutativity of the Hopf algebra, a~much simpler strong
connection formula (which does not need to be iterated, nor requires putting the multipullback in the iterated form~--
the latter being complicated and error prone by itself).
While the assumption of co-commutativity limits severely the applicability of the formula, it is worth pointing out that
many of the known piecewise principal comodule algebras, such as those considered in~\cite{bhms,hms2,hms,Rennie,
haj-wagn,reduct} are either $C({\mathbb Z}_n)$ or $\mathcal{O}(U(1))$-comodule algebras, hence our result
could have been used to compute strong connections for these examples.
The strong connection formula presented in this paper was inspired (very loosely) by the proof
of~\cite[Theorem~3.3.2]{rudnikphd}.

The plan of the paper is as follows: Section~\ref{preliminaries} contains some preliminaries about principal comodule
algebras and piecewise principality.
In Section~\ref{strong-section} we present the explicit formula for a~strong connection, and prove that it is indeed
a~strong connection, as long as the Hopf algebra is co-commutative.
Because the strong connection formula uses the colinear and unital splittings of projections onto pieces, we devote
Section~\ref{colinear-section} to the presentation of the explicit procedure for constructing such splittings from the
appropriate splittings of the gluing maps.
Note that Theorem~\ref{main-split} can be viewed as the strengthening of~\cite[Proposition~9]{CM2000}
(cf.~\cite[Theorem~7]{cocycle})~-- instead of merely showing that, for each element in the multipullback component,
there exists an element in the multipullback projected to this element we explicitly construct the whole (co-)linear and
unital splitting.

As some of the splittings of gluing maps used in the construction of the splitting from Theo\-rem~\ref{main-split} are
required to have fairly non-obvious properties,  Section~\ref{partitions-s} is devoted to showing when such
a~splittings are guaranteed to exist, as well as to their semi-explicit constructions.
Lemma~\ref{partlem}, which links the existence of certain partitions of a~vector space generated by a~collection of
vector subspaces to the distributivity of the lattice generated by those subspaces, is crucial for the results in this
section.

Finally, in Section~\ref{example-section}, we derive a~formula for a~strong connection on a~non-commutative sphere
$S^2_{\mathbb{R}\mathcal{T}}$ introduced in~\cite{reduct} as a~quantum ${\mathbb Z_2}$-principal bundle.
To this end, and to provide comparison, we use two methods~-- the one from~\cite{hkmz11} and the one introduced in this
paper.

\section{Preliminaries}\label{preliminaries}

\subsection{Hopf algebra and comodule-related notation}

We work over a~f\/ixed ground f\/ield $\mathbb{K}$ and, unless stated otherwise, all vector spaces are understood to be
$\mathbb{K}$-vector spaces and the unadorned tensor product is understood to be the algebraic tensor product over~$\mathbb{K}$.
The comultiplication, counit and the antipode of a~Hopf algebra~$H$ are denoted by~$\Delta$, ${\epsilon}$ and~$S$,
respectively.
Let~$P$ be a~right comodule algebra.
We denote by $\Delta_P:P\rightarrow P\otimes H$ the right $H$-coaction on~$P$, and by
\begin{gather*}
P^{\co H}:=\big\{p\in P
\,|\,
\Delta_P(p)=p\otimes 1_H\big\}
\end{gather*}
the subalgebra of coaction invariant elements.
Instead of writing~$\Delta$'s and $\Delta_P$'s we usually employ the Heynemann--Sweedler notation with the summation
symbol suppressed, e.g.,
\begin{gather*}
\Delta(h)=:h\sw{1}\otimes h\sw{2},
\qquad
\Delta_P(p)=:p\sw{0}\otimes p\sw{1}.
\end{gather*}

\subsection{Principal comodule algebras}

Let~$H$ be a~Hopf algebra with bijective antipode, and let~$P$ be a~right $H$-comodule algebra.
Then~$P$ is a~principal comodule algebra if and only if there exists a~linear map
\begin{gather*}
\ell: \ H\rightarrow P\otimes P,
\qquad
\ell(h)=:\ls{h}\otimes\rs{h}
\end{gather*}
(note the Sweedler-like notation with summation sign supressed) satisfying the following conditions
\begin{subequations}
\label{el-prop}
\begin{gather}
\ell(1_H)=1_P\otimes1_P,
\label{el-unitality}
\\
\ls{h}\rs{h}={\epsilon}(h),
\label{el-collapse}
\\
\ls{h\sw{1}}\otimes\rs{h\sw{1}}\otimes h\sw{2}=\ls{h}\otimes\rs{h}\sw{0}\otimes\rs{h}\sw{1},
\label{el-rcolinear}
\\
S(h\sw{1})\otimes\ls{h\sw{2}}\otimes\rs{h\sw{2}}=\ls{h}\sw{1}\otimes\ls{h}\sw{0}\otimes\rs{h}.
\label{el-lcolinear}
\end{gather}
\end{subequations}

Such a~map, if it exists, is called a~{\em strong connection} on~$P$~\cite{bh04,dgh01, h-pm96}.
Strong connections are usually non-unique.

\subsection{Multi-pullbacks of algebras}

Let~$J$ be a~f\/inite set, and let
\begin{gather}
\label{family}
\big\{\pi^i_j:A_i\longrightarrow A_{ij}=A_{ji}\big\}_{i,j\in J,\,i\neq j}
\end{gather}
be a~family of algebra homomorphisms to which we will occasionally refer as ``gluing maps''.
\begin{definition}
[\cite{CM2000,GK-P99}] The \emph{multi-pullback algebra} $A^\pi$ of a~family~\eqref{family} of algebra homomorphisms is
def\/ined as
\begin{gather*}
A^\pi:=\bigg\{ (a_i)_{i\in J}\in\prod\limits_{i\in J}A_i
\,
\bigg|
\,
\pi^i_j(a_i)=\pi^j_i(a_j),
\;
\forall\, i,j\in J,\, i\neq j \bigg\}.
\end{gather*}
\end{definition}

\begin{definition}[\cite{cocycle}] A~family~\eqref{family} of algebra homomorphisms is called \emph{distributive} if and only if all of
them are surjective and their kernels generate distributive lattices of ideals.
\end{definition}

Let $(\pi^i_j:A_i\rightarrow A_{ij})_{i,j\in J,\,i\neq j}$ be a~family of surjective algebra homomorphisms.
For any distinct $i$, $j$, $k$ we put $A^i_{jk}:=A_i/(\ker\pi^i_j+\ker\pi^i_k)$ and take $[\cdot]^i_{jk}:A_i\rightarrow
A^i_{jk}$ to be the canonical surjections.
Next, we introduce the family of maps
\begin{gather*}
\pi^{ij}_k: \ A^i_{jk}\longrightarrow A_{ij}/\pi^i_j\big(\ker\pi^i_k\big),
\qquad
[a_i]^i_{jk}\longmapsto\pi^i_j(a_i)+\pi^i_j\big(\ker\pi^i_k\big).
\end{gather*}
They are isomorphisms when $\pi^i_j$'s are epimorphisms.

\begin{definition}
\label{cocycle-def}
We say \cite[Proposition~9]{CM2000} that a~family $(\pi^i_j:A_i\rightarrow A_{ij})_{i,j\in J,\,i\neq j}$ of algebra
epimorphisms satisf\/ies the {\em cocycle condition} if and only if, for all distinct $i,j,k\in J$,
\begin{enumerate}\itemsep=0pt
\item[1)] $\pi^i_j(\ker\pi^i_k)=\pi^j_i\big(\ker\pi^j_k\big)$,
\item[2)] the isomorphisms $\phi^{ij}_k:=\big(\pi^{ij}_k\big)^{-1}\circ\pi^{ji}_k:A^j_{ik}\rightarrow A^i_{jk}$ satisfy
$\phi^{ik}_j=\phi^{ij}_k\circ\phi^{jk}_i$.
\end{enumerate}
\end{definition}

Observe that, for all distinct $i,j,k\in J$ and any $a_i\in A_i$, $a_j\in A_j$,
\begin{gather}
\label{eqphi}
[a_i]^i_{jk}=\phi^{ij}_k\big([a_j]^j_{ik}\big)
\quad\!
\Leftrightarrow
\quad\!
\pi^{ji}_k\big([a_j]^j_{ik}\big)=\pi^{ij}_k\big([a_i]^{i}_{jk}\big)
\quad\!
\Leftrightarrow
\quad\!
\pi^i_j(a_i)-\pi^j_i(a_j)\in\pi^i_j\big(\ker\pi^i_k\big).\!\!\!
\end{gather}
One can prove (\cite{CM2000}, cf.~\cite{cocycle}, see also Theorem~\ref{main-split} in this paper) that the cocycle
condition together with distributivity guarantees that all projections on components of a~multipullback are surjective
(in fact all projections on submultipullbacks are surjective, but we will not make use of that fact).

\subsection{Piecewise principal comodule algebras}

\begin{definition}[\protect{cf.~\cite[Definition~3.7]{hkmz11}}] A~family of surjective algebra homomorphisms  $\{\pi_i:P\rightarrow P_i\}_{i\in\{1,\ldots,N\}}$
is called a~{\em covering}~\cite{hkmz11} if and only if
\begin{enumerate}\itemsep=0pt
\item[1)] $\bigcap_{i\in\{1,\ldots,N\}}\ker\pi_i=\{0\}$,
\item[2)] the family of ideals $(\ker\pi_i)_{i\in\{1,\ldots,N\}}$ generates a~distributive lattice with $+$ and $\cap$ as
meet and join, respectively.
\end{enumerate}
\end{definition}

Piecewise principal comodule algebras generalize the notion of (algebras of functions on) classical spaces which are
locally principal, but with respect to closed instead of open coverings~-- hence the use of the term ``piecewise''
instead of ``locally''.

\begin{definition}[\protect{see~\cite[Definition~3.8]{hkmz11}}]\label{piecewise-principal-def}
An $H$-comodule algebra~$P$ is called piecewise principal if there exists a~f\/inite
family $\{\pi_i:P\rightarrow P_i\}_{i\in J}$ of surjective $H$-comodule algebra morphisms such that
\begin{enumerate}\itemsep=0pt
\item[1)] the restrictions $\pi_i\big|_{P^{\co H}}:P^{\co H}\rightarrow P^{\co H}_i$ form a~covering,
\item[2)] the $P_i$'s are principal $H$-comodule algebras.
\end{enumerate}
\end{definition}

Note that, for all $i\in J$, $\pi_i(P^{\co H})\subseteq P^{\co H}_i$ by virtue of right $H$-colinearity of~$\pi_i$.
Hence, we were allowed to consider $\pi_i\big|_{P^{\co H}}$ in the statement of
Def\/inition~\ref{piecewise-principal-def} as a~map with codomain~$P^{\co H}_i$ without any additional assumptions.

By~\cite[Corollary~3.9]{hkmz11} a~piecewise principal comodule algebra is principal.
Note that any piecewise principal comodule algebra can be presented as a~multipullback comodule algebra with the gluing
maps being comodule algebra morphisms~\cite{cm2002}.

\section{Strong connection formula}\label{strong-section}

In this section we present an explicit (and arguably simple) expression for a~strong connection on a~piecewise
principal $H$-comodule algebra where~$H$ is a~co-commutative Hopf algebra.
Regretfully, the co-commutativity assumption is used crucially in the proof of the correctness of the formula, and so we
have little hopes of generalizing further the method which led to the derivation of this strong connection formula.

\begin{theorem}
\label{Main}
Let~$H$ be a~cocomutative Hopf algebra.
Let $\{\pi_i:P\rightarrow P_i\}_{i\in\{0,\ldots,n\}}$ be a~piecewise principal $H$-comodule algebra, and let
$\{\ell_i:H\rightarrow P_i\otimes P_i\}_{i\in\{0,\ldots,n\}}$ denote a~family of strong connections on $P_i$'s.
For any $i\in\{0,\ldots,n\}$, let $V_i$ be an~$H$ sub-comodule of $P_i$ such that $\ell_i(H)\subseteq V_i\otimes V_i$
and let $\alpha_i:V_i\rightarrow P$ be a~unital, colinear splitting of $\pi_i$, i.e.,
$\pi_i\circ\alpha_i=\mathrm{id}_{V_i}$.
For brevity, denote for $i\in\{0,\ldots,n\}$, $h\in H$
\begin{gather*}
\theta_i(h):={\epsilon}(h)-\alpha_i\big(\lls{i}{h}\big)\alpha_{i}\big(\rrs{i}{h}\big),
\\
T_i(h):=\theta_i(h\sw{1})\theta_{i+1}(h\sw{2})\cdots\theta_n(h\sw{n-i+1}),
\qquad
T_{n+1}(h):={\epsilon}(h).
\end{gather*}
Then the linear map $\ell:H\rightarrow P\otimes P$ defined for all $h\in H$ by the formula
\begin{gather*}
\ell(h)=\sum\limits_{i=0}^n\alpha_i\big(\lls{i}{h\sw{1}}\big)\otimes\alpha_i\big(\rrs{i}{h\sw{1}}\big)T_{i+1}(h\sw{2})
\end{gather*}
is a~strong connection on~$P$.
\end{theorem}

Note that, in particular, $T_n(h)=\theta_n(h)$, for all $h\in H$.
Note also that we consider splittings from $V_i$'s instead of splittings from $P_i$'s because the former are much easier
to construct.
\begin{proof}
Note that any co-commutative Hopf algebra has bijective (in fact involutive) antipode.
We need to prove that the map~$\ell$ def\/ined in the theorem, satisf\/ies all the properties~\eqref{el-prop}.

First note that, by the colinearity of $\alpha_j$'s, colinear properties~\eqref{el-lcolinear},~\eqref{el-rcolinear} of
$\ell_j$'s and the co-commutativity of~$H$ we have, that $\alpha_j(\lls{j}{h})\alpha_{j}(\rrs{j}{h})$ is a~coaction
invariant element of~$P$ for any $j\in\{0,\ldots,n\}$ and $h\in H$, and hence also $T_i(h)$ is a~coaction invariant
element of~$P$ for any $i\in\{0,\ldots,n+1\}$ and $h\in H$
\begin{gather*}
\rho^H\big(\alpha_j\big(\lls{j}{h}\big)\alpha_{j}\big(\rrs{j}{h}\big)\big)\\
\qquad{}
=\alpha_j\big(\lls{j}{h}\big)\sw{0}\alpha_{j}\big(\rrs{j}{h}\big)\sw{0} \otimes
\alpha_j\big(\lls{j}{h}\big)\sw{1}\alpha_{j}\big(\rrs{j}{h}\big)\sw{1}
\\
\qquad{}
 =\alpha_j\big(\lls{j}{h\sw{2}}\big)\alpha_{j}\big(\rrs{j}{h\sw{2}}\big)\otimes S(h\sw{1})h\sw{3}
\\
\qquad{}
 =\alpha_j\big(\lls{j}{h\sw{1}}\big)\alpha_{j}\big(\rrs{j}{h\sw{1}}\big)\otimes S(h\sw{2})h\sw{3}\\
 \qquad{}
 =\alpha_j\big(\lls{j}{h}\big)\alpha_{j}\big(\rrs{j}{h}\big)\otimes1.
\end{gather*}
In the penultimate equality we used co-commutativity of~$H$ to swap Sweedler indices $\sw{1}$ and $\sw{2}$ to be able to
use the antipode property.
In order to prove that~$\ell$ is left colinear (equation~\eqref{el-lcolinear}) we use the left colinearity of $\ell_i$'s
and the right colinearity of $\alpha_i$'s
\begin{gather*}
\ls{h}\sw{1}\otimes\ls{h}\sw{0}\otimes\rs{h}
\\
\qquad=
\sum\limits_{i=0}^n\alpha_i\big(\lls{i}{h\sw{1}}\big)\sw{1}\otimes\alpha_i\big(\lls{i}{h\sw{1}}\big)\sw{0}
\otimes\alpha_i\big(\rrs{i}{h\sw{1}}\big)T_{i+1}(h\sw{2})
\\
\qquad=
\sum\limits_{i=0}^n\lls{i}{h\sw{1}}\sw{1}\otimes\alpha_i(\lls{i}{h\sw{1}}\sw{0})
\otimes\alpha_i\big(\rrs{i}{h\sw{1}}\big)T_{i+1}(h\sw{2})
\\
\qquad=
\sum\limits_{i=0}^nS(h\sw{1})\otimes\alpha_i\big(\lls{i}{h\sw{2}}\big) \otimes\alpha_i\big(\rrs{i}{h\sw{2}}\big)T_{i+1}(h\sw{3})
\\
\qquad=
S(h\sw{1})\otimes\ls{h\sw{2}}\otimes\rs{h\sw{2}}.
\end{gather*}
The right colinearity 
(equation~\eqref{el-rcolinear}) of~$\ell$ follows from the $H$-coaction invariance of $T_i(h)$'s,
the right colinearity of $\ell_i$'s, the right colinearity of $\alpha_i$'s, and the co-commutativity of~$H$
\begin{gather*}
\ls{h}\otimes\rs{h}\sw{0}\otimes\rs{h}\sw{1}
\\
\qquad=
\sum\limits_{i=0}^n\alpha_i\big(\lls{i}{h\sw{1}}\big)\otimes
\alpha_i\big(\rrs{i}{h\sw{1}}\big)\sw{0}T_{i+1}(h\sw{2})\otimes\alpha_i\big(\rrs{i}{h\sw{1}}\big)\sw{1}
\\
\qquad=
\sum\limits_{i=0}^n\alpha_i\big(\lls{i}{h\sw{1}}\big)\otimes
\alpha_i\big(\rrs{i}{h\sw{1}}\sw{0}\big)T_{i+1}(h\sw{2})\otimes\rrs{i}{h\sw{1}}\sw{1}
\\
\qquad=
\sum\limits_{i=0}^n\alpha_i\big(\lls{i}{h\sw{1}}\big)\otimes \alpha_i\big(\rrs{i}{h\sw{1}}\big)T_{i+1}(h\sw{3})\otimes h\sw{2}
\\
\qquad=
\sum\limits_{i=0}^n\alpha_i\big(\lls{i}{h\sw{1}}\big)\otimes \alpha_i\big(\rrs{i}{h\sw{1}}\big)T_{i+1}(h\sw{2})\otimes h\sw{3}
\\
\qquad=
\ls{h\sw{1}}\otimes\rs{h\sw{1}}\otimes h\sw{3}.
\end{gather*}
Here, in the penultimate inequality we used the co-commutativity of~$H$ exchanging Sweedler indices $\sw{2}$ and $\sw{3}$.

In order to prove that~$\ell$ is unital (equation~\eqref{el-unitality}), note f\/irst that for any $i\in\{0,\ldots,n\}$
\begin{gather*}
\theta_i(1)={\epsilon}(1)-\alpha_i\big(\lls{i}{1}\big)\alpha_{i}\big(\rrs{i}{1}\big)=1-\alpha_i(1)\alpha_i(1)=1-1=0,
\end{gather*}
because ${\epsilon}$, all $\ell_i$'s and all $\alpha_i$'s are unital.
It follows that $T_i(1)=0$ for all $i\in\{0,\ldots,n\}$, and $T_{n+1}={\epsilon}$ by def\/inition, hence
\begin{gather*}
\ell(1)=\sum\limits_{i=0}^n\!\alpha_i\big(\lls{i}{1}\big)\otimes\alpha_i\big(\rrs{i}{1}\big)T_{i+1}(1)
=\alpha_n\big(\lls{n}{1}\big)\otimes\alpha_n\big(\rrs{n}{1}\big)T_{n+1}(1)=1\otimes1,
\end{gather*}
where we used again the unitality of $\alpha_n$ and $\ell_n$.

Note now that for all $i\in\{0,\ldots,n\}$, and $h\in H$
\begin{gather*}
T_i(h)=T_{i+1}(h)-\alpha_i\big(\lls{i}{h\sw{1}}\big)\alpha_i\big(\rrs{i}{h\sw{1}}\big)T_{i+1}(h\sw{2}).
\end{gather*}
Indeed,
\begin{gather*}
T_i(h) =\theta_i(h\sw{1})T_{i+1}(h\sw{2})
 ={\epsilon}(h\sw{1})T_{i+1}(h\sw{2})-\alpha_i\big(\lls{i}{h\sw{1}}\big)\alpha_i\big(\rrs{i}{h\sw{1}}\big)T_{i+1}(h\sw{2})
\\
\phantom{T_i(h)}
 =T_{i+1}(h)-\alpha_i\big(\lls{i}{h\sw{1}}\big)\alpha_i\big(\rrs{i}{h\sw{1}}\big)T_{i+1}(h\sw{2}).
\end{gather*}
By applying this formula to $T_0(h)$ and keeping to expand with it the leftmost summand of the resulting expansion we
obtain easily
\begin{gather}
\label{t0}
T_0(h)={\epsilon}(h)-\sum\limits_{i=0}^n\alpha_i\big(\lls{i}{h\sw{1}}\big)\alpha_i\big(\rrs{i}{h\sw{1}}\big)T_{i+1}(h\sw{2}).
\end{gather}
On the other hand, for all $h\in H$ and $i\in\{0,\ldots,n\}$, as $\alpha_i$ is the splitting of $\pi_i$ it follows that
\begin{gather*}
\pi_i(\theta_i(h))  ={\epsilon}(h)-\pi_i\big(\alpha_i(\lls{i}{h})\big)\pi_i\big(\alpha_i(\rrs{i}{h})\big)
\\
\phantom{\pi_i(\theta_i(h))}
={\epsilon}(h)-\lls{i}{h}\rrs{i}{h} ={\epsilon}(h)-{\epsilon}(h)
 =0.
\end{gather*}
Hence
\begin{gather*}
\pi_i(T_j(h))=0,
\qquad
\text{for all}
\quad
i\geq j,
\quad
i\in\{0,\ldots,n\},
\quad
h\in H.
\end{gather*}
In particular, $\pi_i(T_0(h))=0$ for all $i\in\{0,\ldots,n\}$ and $h\in H$.
It follows that $T_0(h)=0$ for all $h\in H$ because $\bigcap_{i=0}^n\ker\pi_i=\{0\}$, as $\{\pi_i:P\rightarrow
P_i\}_{i\in\{0,\ldots,n\}}$ is a~covering.
The last fact is an immediate consequence of~\cite[Theorem~3.3]{hkmz11} and~\cite[Corollary~3.7]{hkmz11}.

Combining this with the equation~\eqref{t0} we obtain that for all $h\in H$
\begin{gather*}
\ls{h}\rs{h}=\sum\limits_{i=0}^n\alpha_i\big(\lls{i}{h\sw{1}}\big)\alpha_i\big(\rrs{i}{h\sw{1}}\big)T_{i+1}(h\sw{2})={\epsilon}(h),
\end{gather*}
i.e.,~$\ell$ satisf\/ies equation~\eqref{el-collapse} as needed.
\end{proof}

The expression for a~strong connection provided in the above theorem requires the unital and colinear splittings of
projections $\pi_i$ to be given.
The existence of such a~splittings is guaranteed by the~\cite[Lemma~3.1]{hkmz11} and~\cite[Theorem~3.3]{hkmz11}, but the
mere existence does not suf\/f\/ice for someone desirous of f\/inding the explicit formula.
The proof of~\cite[Lemma~3.1]{hkmz11} involves constructing a~unital and colinear splitting of surjective comodule
algebra map~$\pi$ from a~unital and linear splitting of restriction of~$\pi$ to the subalgebra of coaction invariant
elements (which always exists) utilizing the strong connection.
Hence, we cannot use even the slight simplif\/ication provided by the proof of~\cite[Lemma~3.1]{hkmz11}.

In practice, we expect that in many simpler cases, the appropriate splittings will not be dif\/f\/icult to guess.
However, for our result to be more widely applicable in practice, we will examine the explicit construction of colinear
and unital splittings of multipullback comodule algebra projections on components which does not assume the existence of
a~strong connection on a~multipullback comodule algebra (recall that a~piecewise principal comodule algebra can always
be presented as a~multipullback).

\section{Colinear splittings of piecewise principal comodule algebras}\label{colinear-section}

The result presented in this section allows to explicitly construct linear (colinear when approp\-riate) and unital
splittings of projections on components of a~multipullback (comodule) algebra.

\begin{theorem}
\label{main-split}
Suppose that a~family~\eqref{family} is distributive and satisfies the cocycle condition.
Moreover suppose that there exists two families $\alpha^i_j,\beta^i_j:A_{ij}\rightarrow A_i$, $i,j\in J$,
$j\neq i$ of linear $($colinear$)$ splittings of $\pi^i_j$'s such that all $\beta^i_j$'s are unital and for all distinct
$i,j,k\in J$ we have
\begin{gather}
\label{thkercond}
\alpha^i_j\big(\pi^i_j\big(\ker\pi^i_k\big)\big)\subseteq\ker\pi^i_k.
\end{gather}
Let $i\in J$, let $|J|=n+1$ and let $\kappa:\{0,\ldots,n\}\rightarrow J$ be a~bijection such that $\kappa_0=i$, where we
denote $\kappa_j:=\kappa(j)$ to easy the notation.
Then a~unital and linear $($colinear$)$ splitting $\alpha_i:A_i\rightarrow A^\pi$ of $\pi_i:A^\pi\rightarrow A_i$ can be
given explicitly, for any $a\in A_i$ as $\alpha_i(a):=(a_j)_{j\in J}$, where $a_i:=a$ and
$a_{\kappa_{m+1}}:=a^m_{\kappa_{m+1}}$ for any $0\leq m<n$.
The collections $\{a^k_{\kappa_{m+1}}\}_{0\leq k\leq m}\subseteq A_{\kappa_{m+1}}$, for $0\leq m<n$ are defined by the
following inductive formula
\begin{gather}
a^0_{\kappa_{m+1}}:=\beta^{\kappa_{m+1}}_{\kappa_0}\big(\pi^{\kappa_0}_{\kappa_{m+1}}(a_{\kappa_0})\big),
\nonumber
\\
a^{k+1}_{\kappa_{m+1}}:=a_{\kappa_{m+1}}^k-\alpha^{\kappa_{m+1}}_{\kappa_{k+1}}
\big(\pi^{\kappa_{m+1}}_{\kappa_{k+1}}\big(a^k_{\kappa_{m+1}}\big)-\pi^{\kappa_{k+1}}_{\kappa_{m+1}}(a_{\kappa_{k+1}})\big)
\label{recsplitform}
\end{gather}
for $0\leq k<m$.
\end{theorem}

\begin{proof}
It is clear that because all the maps involved in the def\/inition of $\alpha_i$ are unital and linear (colinear if need
be) then also $\alpha_i$ is linear (resp.
colinear).
The proof of unitality is slightly more subtle and it requires a~simple induction.
Pick some bijection $\kappa:\{0,\ldots,n\}\rightarrow J$ where $\kappa_0=i$.
Def\/ine $(a_j)_{j\in J}:=\alpha_i(1)$.
We need to show that $a_j=1$ for all $j\in J$.
Indeed, $a_{\kappa_0}=a_i=1$ by def\/inition.
Suppose we have proven that $a_j=1$ for all $0\leq j\leq m<n$.
Then using the equation~\eqref{recsplitform} we get
$a^0_{\kappa_{m+1}}=\beta^{\kappa_{m+1}}_{\kappa_0}(\pi^{\kappa_0}_{\kappa_{m+1}}(a_{\kappa_0}))=
\beta^{\kappa_{m+1}}_{\kappa_0}(\pi^{\kappa_0}_{\kappa_{m+1}}(1))=1$ as both~$\pi^{\kappa_0}_{\kappa_{m+1}}$ and~$\beta^{\kappa_{m+1}}_{\kappa_0}$ are unital.
Suppose now that we have proven that $a^k_{\kappa_{m+1}}=1$ for all $0\leq k<m$.
Then, equation~\eqref{recsplitform} yields
\begin{gather*}
a^{k+1}_{\kappa_{m+1}}
=a_{\kappa_{m+1}}^k-\alpha^{\kappa_{m+1}}_{\kappa_{k+1}}
\big(\pi^{\kappa_{m+1}}_{\kappa_{k+1}}(a^k_{\kappa_{m+1}})-\pi^{\kappa_{k+1}}_{\kappa_{m+1}}(a_{\kappa_{k+1}})\big)
\\
\phantom{a^{k+1}_{\kappa_{m+1}}}
 =1 -\alpha^{\kappa_{m+1}}_{\kappa_{k+1}}\big(\pi^{\kappa_{m+1}}_{\kappa_{k+1}}(1)-\pi^{\kappa_{k+1}}_{\kappa_{m+1}}(1)\big)
 =1-\alpha^{\kappa_{m+1}}_{\kappa_{k+1}}(0)
 =1.
\end{gather*}

Now it remains to show that $\alpha_i(a)\in A^\pi$ for all $a\in A_i$.
The inductive proof essentially follows the steps of the proof of~\cite[Proposition~9]{CM2000}.
We will show that for any $0\leq m\leq n$ we have
\begin{gather}
\label{partcond1}
\pi^{\kappa_j}_{\kappa_l}(a_{\kappa_j})=\pi^{\kappa_l}_{\kappa_j}(a_{\kappa_l}),
\qquad
\text{for all}
\quad
j,l\in\{0,\ldots,m\},
\quad
j\neq l.
\end{gather}
For $m=0$ this condition is emptily satisf\/ied.
Suppose we have proven the above condition for some~$m$.
In order to demonstrate it for $m+1$, we prove by induction that for any $0\leq k\leq m$, where $m<n$, we have
\begin{gather}
\label{partcond2}
\pi^{\kappa_j}_{\kappa_{m+1}}(a_{\kappa_j})=\pi^{\kappa_{m+1}}_{\kappa_j}\big(a^k_{\kappa_{m+1}}\big),
\qquad
\text{for all}
\quad
0\leq j\leq k.
\end{gather}
If $k=0$ then substituting the def\/inition of $a^0_{\kappa_{m+1}}$ yields (as $\beta^{\kappa_{m+1}}_{\kappa_0}$ is
a~splitting of $\pi^{\kappa_{m+1}}_{\kappa_0}$)
\begin{gather*}
\pi^{\kappa_{m+1}}_{\kappa_0}(a^0_{\kappa_{m+1}}) =\pi^{\kappa_{m+1}}_{\kappa_0}
\big(\beta^{\kappa_{m+1}}_{\kappa_0}\big(\pi^{\kappa_0}_{\kappa_{m+1}}(a_{\kappa_0})\big)\big)
 =\pi^{\kappa_0}_{\kappa_{m+1}}(a_{\kappa_0}).
\end{gather*}
Suppose now that we have proven condition~\eqref{partcond2} for some $0\leq k<m$.
Pick any $0\leq j\leq k$.
Then by (inductively assumed) condition~\eqref{partcond1} and equation~\eqref{eqphi} we have
\begin{gather}
\label{usedcoc}
[a_{\kappa_j}]^{\kappa_j}_{\kappa_{k+1}\kappa_{m+1}}=\phi^{\kappa_j\kappa_{k+1}}_{\kappa_{m+1}}
\big([a_{\kappa_{k+1}}]^{\kappa_{k+1}}_{\kappa_j\kappa_{m+1}}\big).
\end{gather}
Then it follows that
\begin{gather*}
\big[a^k_{\kappa_{m+1}}\big]^{\kappa_{m+1}}_{\kappa_j\kappa_{k+1}}
\overset{\substack{\text{by condition~\eqref{partcond2}} \\ \text{and equation~\eqref{eqphi}}}}{=}\phi^{\kappa_{m+1}\kappa_j}_{\kappa_{k+1}}
\big([a_{\kappa_j}]^{\kappa_j}_{\kappa_{m+1}\kappa_{k+1}}\big)
\\
\qquad{}
\overset{\text{by equation~\eqref{usedcoc}}}{=}
\phi^{\kappa_{m+1}\kappa_j}_{\kappa_{k+1}}\big(\phi^{\kappa_j\kappa_{k+1}}_{\kappa_{m+1}}
 \big([a_{\kappa_{k+1}}]^{\kappa_{k+1}}_{\kappa_j\kappa_{m+1}}
\big) \big)
\overset{\substack{\text{by the cocycle} \\ \text{condition\vphantom{q}}}}{=} \phi^{\kappa_{m+1}\kappa_{k+1}}_{\kappa_j}\big([a_{\kappa_{k+1}}]^{\kappa_{k+1}}_{\kappa_j\kappa_{m+1}}\big).
\end{gather*}
This equality, again by equation~\eqref{eqphi}, is equivalent to the following condition
\begin{gather*}
\pi^{\kappa_{m+1}}_{\kappa_{k+1}}\big(a^k_{\kappa_{m+1}}\big)-\pi^{\kappa_{k+1}}_{\kappa_{m+1}}(a_{\kappa_{k+1}})
\in\pi^{\kappa_{m+1}}_{\kappa_{k+1}}\big(\ker\pi^{\kappa_{m+1}}_{\kappa_j}\big).
\end{gather*}
Because the above relation ``is an element of'' holds for an arbitrary $0\leq j\leq k$ it implies immediately that
\begin{gather}
\label{kercond}
\pi^{\kappa_{m+1}}_{\kappa_{k+1}}\big(a^k_{\kappa_{m+1}}\big)-\pi^{\kappa_{k+1}}_{\kappa_{m+1}}(a_{\kappa_{k+1}})
\in\bigcap_{0\leq j\leq k}\pi^{\kappa_{m+1}}_{\kappa_{k+1}}\big(\ker\pi^{\kappa_{m+1}}_{\kappa_j}\big).
\end{gather}
Then
\begin{gather*}
\alpha^{\kappa_{m+1}}_{\kappa_{k+1}}\big(\pi^{\kappa_{m+1}}_{\kappa_{k+1}}\big(a^k_{\kappa_{m+1}}\big)
-\pi^{\kappa_{k+1}}_{\kappa_{m+1}}(a_{\kappa_{k+1}})\big)
\overset{\text{by condition~\eqref{kercond}}}{\in}
\alpha^{\kappa_{m+1}}_{\kappa_{k+1}}\left(\bigcap_{0\leq j\leq k}\pi^{\kappa_{m+1}}_{\kappa_{k+1}}\big(\ker\pi^{\kappa_{m+1}}_{\kappa_j}\big) \right)
\\
\qquad
\overset{\text{by injectivity of $\alpha^{\kappa_{m+1}}_{\kappa_{k+1}}$}}{\in}
\bigcap_{0\leq j\leq k}
\alpha^{\kappa_{m+1}}_{\kappa_{k+1}}\left(\pi^{\kappa_{m+1}}_{\kappa_{k+1}}\big(\ker\pi^{\kappa_{m+1}}_{\kappa_j}\big)\right)
\overset{\text{by equation~\eqref{thkercond}}}{\subseteq}
\bigcap_{0\leq j\leq k}\ker\pi^{\kappa_{m+1}}_{\kappa_j},
\end{gather*}
that is
\begin{gather*}
\alpha^{\kappa_{m+1}}_{\kappa_{k+1}}\big(\pi^{\kappa_{m+1}}_{\kappa_{k+1}}\big(a^k_{\kappa_{m+1}}\big)
-\pi^{\kappa_{k+1}}_{\kappa_{m+1}}(a_{\kappa_{k+1}})
\big)\in \bigcap_{0\leq j\leq k}\ker\pi^{\kappa_{m+1}}_{\kappa_j}.
\end{gather*}
The above equation implies immediately, that for all $0\leq l\leq k$
\begin{gather*}
\pi^{\kappa_{m+1}}_{\kappa_l}\big(a^{k+1}_{\kappa_{m+1}}\big)
= \pi^{\kappa_{m+1}}_{\kappa_l}\big(a_{\kappa_{m+1}}^k\big)-\pi^{\kappa_{m+1}}_{\kappa_l}\big(\alpha^{\kappa_{m+1}}_{\kappa_{k+1}}
\big(\pi^{\kappa_{m+1}}_{\kappa_{k+1}}\big(a^k_{\kappa_{m+1}}\big)-\pi^{\kappa_{k+1}}_{\kappa_{m+1}}(a_{\kappa_{k+1}})\big)\big)
\\
\phantom{\pi^{\kappa_{m+1}}_{\kappa_l}\big(a^{k+1}_{\kappa_{m+1}}\big)}{}
 =\pi^{\kappa_{m+1}}_{\kappa_l}\big(a_{\kappa_{m+1}}^k\big)
 =\pi^{\kappa_l}_{\kappa_{m+1}}(a_{\kappa_l}),
\end{gather*}
where, in the second equality we used the inductive assumption.
Moreover, using the fact that $\alpha^{\kappa_{m+1}}_{\kappa_{k+1}}$ is a~splitting of
$\pi^{\kappa_{m+1}}_{\kappa_{k+1}}$ we obtain
\begin{gather*}
\pi^{\kappa_{m+1}}_{\kappa_{k+1}}\big(a^{k+1}_{\kappa_{m+1}}\big)
= \pi^{\kappa_{m+1}}_{\kappa_{k+1}}\big(a_{\kappa_{m+1}}^k\big)-\pi^{\kappa_{m+1}}_{\kappa_{k+1}}\big(\alpha^{\kappa_{m+1}}_{\kappa_{k+1}}
\big(\pi^{\kappa_{m+1}}_{\kappa_{k+1}}\big(a^k_{\kappa_{m+1}}\big)-\pi^{\kappa_{k+1}}_{\kappa_{m+1}}(a_{\kappa_{k+1}})\big)\big)
\\
\phantom{\pi^{\kappa_{m+1}}_{\kappa_{k+1}}\big(a^{k+1}_{\kappa_{m+1}}\big)}
 =\pi^{\kappa_{m+1}}_{\kappa_{k+1}}\big(a_{\kappa_{m+1}}^k\big) -
\big(\pi^{\kappa_{m+1}}_{\kappa_{k+1}}\big(a^k_{\kappa_{m+1}}\big)-\pi^{\kappa_{k+1}}_{\kappa_{m+1}}(a_{\kappa_{k+1}})\big)
 =\pi^{\kappa_{k+1}}_{\kappa_{m+1}}(a_{\kappa_{k+1}}),
\end{gather*}
which ends the proof.
\end{proof}

At this point, the skeptical reader might be excused for doubting the applicability of Theorem~\ref{main-split}.
Indeed, while the existence of unital and linear splittings $\beta^i_j$'s of $\pi^i_j$'s follows immediately from the
surjectivity of $\pi^i_j$'s, and the existence of colinear splittings is assured (and assisted in explicit construction)
by~\cite[Lemma~3.1]{hkmz11} if all the $A_i$'s are principal comodule algebras, it is not clear how to f\/ind the linear
splittings $\alpha^i_j$ satisfying equation~\eqref{thkercond} nor that they exist at all in general case.
Fortunately, the results from the next section, interesting in their own right, not only assure the existence of
splittings $\alpha^i_j$ satisfying equation~\eqref{thkercond} under no stronger assumptions than those of
Theorem~\ref{main-split}, but they also provide the method of their (semi)-explicit construction.

\section{Colinear splittings of principal comodule algebras}\label{partitions-s}

\subsection{Partitions of sets}

Let~$A$ be a~set and let $A_i$, $i\in J$ be a~f\/ixed f\/inite family of subsets of~$A$.
For any $\Gamma\in 2^J$ we denote for brevity
\begin{gather}
\label{agamma}
A_\Gamma:=\bigcap_{i\in \Gamma}A_i.
\end{gather}
Obviously $A_{\Gamma_1}\cap A_{\Gamma_2}=A_{\Gamma_1\cup\Gamma_2}$.
Also $A_\varnothing = A$ by convention.
It is easy to see that $A_i$'s generate a~partition $\{B_\Gamma\}_{\Gamma\in 2^J}$ of~$A$ (i.e., all $B_\Gamma$'s are
disjoint and $A=\bigcup_{\Gamma\in2^J}B_\Gamma$) such that
\begin{gather*}
A_\Gamma=\bigcup_{\Gamma'\in2^J
\,| \,
\Gamma\subseteq\Gamma'}B_{\Gamma'},
\qquad
\text{for all}
\quad
\Gamma\in 2^J.
\end{gather*}
Indeed, the partition can be described explicitly, for all $\Gamma\in2^J$ by the formula
\begin{gather*}
B_\Gamma
:=\{x\in A
\,|\,
\forall\, i\in J,\,
x\in A_i\Leftrightarrow i\in \Gamma\}.
\end{gather*}

\subsection{Partitions of vector spaces}

Let now~$A$ be a~vector space and let $A_i$, $i\in J$ be a~f\/ixed f\/inite family of vector subspaces of~$A$.
$A_\Gamma$, for any $\Gamma\in 2^J$ is def\/ined as in equation~\eqref{agamma}.
We want to def\/ine a~linear counterpart of an associated partition $\{B_\Gamma\}_\Gamma$ def\/ined above for sets.
Similarly to plain sets, vector sub-spaces can be ordered by the set inclusion, and the resulting ordered set is
a~lattice, with subspace intersection ($V_1\cap V_2$) serving as inf\/imum and subspace sum ($V_1+V_2$) playing the role
of supremum.
The problem is that this lattice is not, in general, distributive.
It turns out that the assumption that the subspaces $A_i$, $i\in J$ generate a~distributive lattice is pivotal for
proving our desired result, stated immediately below:

\begin{lemma}
\label{partlem}
Let~$A$ be a~linear vector space and let $A_i$, $i\in I$ be a~finite family of vector subspaces of~$A$ generating a~distributive lattice.
$A$ has a~linear basis $\mathcal{B}=\bigcup_{\Gamma\in2^I}\mathcal{B}_\Gamma$, where
$\mathcal{B}_\Gamma\subseteq A_\Gamma$, $\Gamma\in2^I$, such that subsets $\mathcal{B}_\Gamma$ are all disjoint and
satisfy the following property
\begin{gather}
\label{partition}
A_\Gamma=\Span\left(\bigcup_{\Gamma'\in2^I,\; \Gamma'\supseteq\Gamma}\mathcal{B}_{\Gamma'}\right)  
\end{gather}
for all $\Gamma\in2^I$.
\end{lemma}
\begin{proof}
First f\/ix a~linear order $\leq$ on $2^I$ subject to the condition
\begin{gather}
\label{orderdef}
\Gamma_1\supseteq \Gamma_2
\quad
\Rightarrow
\quad
\Gamma_1\leq\Gamma_2,
\qquad
\text{for all}
\quad
\Gamma_1,\Gamma_2\in 2^I.
\end{gather}
It is immediate that the minimal element in this order is~$I$ and maximal is $\varnothing$.
Note the following property of $\leq$ which will be used later
\begin{gather}
\label{Order}
\Gamma>\Gamma'
\quad
\Rightarrow
\quad
\Gamma\cup\Gamma'\supset \Gamma,
\qquad
\text{for all}
\quad
\Gamma,\Gamma'\in2^I.
\end{gather}
Indeed, assume $\Gamma>\Gamma'$.
$\Gamma\cup\Gamma'\supseteq \Gamma$ always, so we need just to show that the equality leads to contradiction.
Suppose that $\Gamma\cup\Gamma'= \Gamma$.
This is equivalent to $\Gamma\supseteq\Gamma'$ which implies by equation~\eqref{orderdef} that $\Gamma\leq\Gamma'$
contradicting the assumption $\Gamma>\Gamma'$.

The sets $\mathcal{B}_\Gamma$, $\Gamma\in2^I$ can be generated inductively (with respect to $\leq$) as follows
\begin{enumerate}\itemsep=0pt
\item[1)] $\mathcal{B}_I$ is some linear basis of $A_I$,
\item[2)] $\mathcal{B}_\Gamma$, for $\Gamma>I$, is chosen as a~maximal subset of $A_\Gamma$ such that
$\bigcup_{\Gamma'\leq\Gamma}\mathcal{B}_{\Gamma'}$ is linearly independent.
\end{enumerate}
It is immediate by construction of $\mathcal{B}_\Gamma$'s that $\mathcal{B}:=\bigcup_{\Gamma\in2^I}\mathcal{B}_\Gamma$
is a~linear basis of~$A$ and that all $\mathcal{B}_\Gamma$'s are disjoint.
Also by construction, $\mathcal{B}_{\Gamma'}\subseteq A_\Gamma$, $\Gamma\in2^I$ whenever $\Gamma\subseteq \Gamma'$,
which implies that half of property~\eqref{partition} is trivially satisf\/ied:
\begin{gather*}
\Span\left(\bigcup_{\Gamma'\in2^I,\, \Gamma'\supseteq\Gamma}\mathcal{B}_{\Gamma'}\right)\subseteq A_\Gamma
\end{gather*}
for all $\Gamma\in2^I$.
Finally, it is immediate that
\begin{gather}
\label{partition0}
A_\Gamma\subseteq\Span\left(\bigcup_{\Gamma'\in2^I,\, \Gamma'\leq\Gamma}\mathcal{B}_{\Gamma'}\right).
\end{gather}
We will prove the second half of property~\eqref{partition} by induction on $\leq$.

1.~$I$ is minimal in $2^I$ with respect to $\leq$.
Then by def\/inition of $B_I$ we have
\begin{gather*}
A_I= \Span(\mathcal{B}_I) = \Span\left(\bigcup_{\Gamma'\in2^I,\, \Gamma'\supseteq I}\mathcal{B}_{\Gamma'}\right).
\end{gather*}

2.~Suppose we have proven equation~\eqref{partition} for all $\Gamma'<\Gamma$.
For any $a\in A$, denote by $\{\alpha_\Gamma(a)\}_{\Gamma\in 2^I}$ the unique family of vectors such that
$a=\sum\limits_{\Gamma\in2^I}\alpha_\Gamma(a)$ and that $\alpha_\Gamma(a)\in \Span(\mathcal{B}_\Gamma)$ for all
$\Gamma\in2^I$ (they are unique because $\mathcal{B}$ is a~basis and $\mathcal{B}_\Gamma$'s are disjoint).
By~\eqref{partition0} $\alpha_{\Gamma'}(a)=0$ whenever $a\in A_\Gamma$ and $\Gamma'>\Gamma$, i.e.,
\begin{gather}
\label{partbase}
a=\sum\limits_{\Gamma'\in2^I,\, \Gamma'\leq\Gamma}\alpha_{\Gamma'}(a),
\qquad
\text{for all}
\quad
a\in A_\Gamma.
\end{gather}
Let $a\in A_\Gamma$.
Def\/ine $v:=a-\alpha_\Gamma(a)$.
By equation~\eqref{partbase}
\begin{gather*}
A_\Gamma\ni v = \sum\limits_{\Gamma'\in 2^I,
\,
\Gamma'<\Gamma}\alpha_{\Gamma'}(a)\in\sum\limits_{\Gamma'\in 2^I,
\,
\Gamma'<\Gamma}A_{\Gamma'},
\end{gather*}
hence
\begin{gather*}
v \in A_\Gamma\cap\left(\sum\limits_{\Gamma'\in 2^I,
\,
\Gamma'<\Gamma}A_{\Gamma'}\right)
\overset{\substack{\text{by distributivity of lattice}\\ \text{generated by $A_i$'s}}}{=}
\sum\limits_{\Gamma'\in 2^I,
\,
\Gamma'<\Gamma}A_{\Gamma'\cup\Gamma}
\\
\qquad \overset{\substack{\text{by equation~\eqref{Order}} \\ \text{$\Gamma\subset\Gamma\cup\Gamma'$ if $\Gamma'<\Gamma$}}}{\subseteq}
\sum\limits_{\Gamma'\in 2^I,
\,
\Gamma'\supset\Gamma}A_{\Gamma'}
 \overset{\substack{\text{by inductive assumption,} \\ \text{as $\Gamma'<\Gamma$ if $\Gamma'\supset\Gamma$}}}{\subseteq}
\Span\left(\bigcup_{\Gamma'\in 2^I,
\,
\Gamma'\supset\Gamma}\mathcal{B}_{\Gamma'} \right).
\end{gather*}
It follows that
\begin{gather*}
a= \alpha_\Gamma(a)+v\in\Span(\mathcal{B}_\Gamma)+\Span\left(\bigcup_{\Gamma'\in 2^I,
\,
\Gamma'\supset\Gamma}\mathcal{B}_{\Gamma'} \right) = \Span\left(\bigcup_{\Gamma'\in 2^I,
\,
\Gamma'\supseteq\Gamma}\mathcal{B}_{\Gamma'} \right)
\end{gather*}
as needed.
\end{proof}

The following result is a~common knowledge:
\begin{lemma}
\label{common}
Let $\pi:A\rightarrow B$ be a~linear map, and let $\{A_i\}_{i\in I}$ be a~finite family of vector subspaces of~$A$.
Assume that $\ker\pi\cap\Big(\sum\limits_{i\in I}A_i\Big)=\sum\limits_{i\in I}(\ker\pi\cap A_i)$.
Then
\begin{gather*}
\pi\bigg(\bigcap_{i\in I}A_i\bigg)=\bigcap_{i\in I}\pi(A_i).
\end{gather*}
\end{lemma}

\begin{lemma}
\label{splitcoinv}
Let $\pi:A\rightarrow B$ be a~linear epimorphism, and let $\{A_i\}_{i\in I}$ be a~finite family of vector subspaces
of~$A$ such that $\{A_i\}_{i\in I}\cup\{\ker\pi\}$ generates a~distributive lattice of vector subspaces.
Then there exists a~linear splitting $\alpha:B\rightarrow A$ of~$\pi$ such that $\alpha(\pi(A_i))\subseteq A_i$ for all
$i\in I$.
\end{lemma}
\begin{proof}
Let $\mathcal{B}:=\bigcup_{\Gamma\in 2^I}\mathcal{B}_\Gamma$ be a~linear basis of~$B$ satisfying conditions guaranteed
by Lemma~\ref{partlem} with respect to the family $\{B_i\}_{i\in I}$, where $B_i:=\pi(A_i)$.
Note that Lemma~\ref{common} implies that $B_i$'s generate distributive lattice of ideals because $A_i$'s generate
distributive lattice of ideals.
For all $\Gamma\in2^I$ such that $\mathcal{B}_\Gamma$ is non-empty we def\/ine $\alpha(b)$ for all $b\in
\mathcal{B}_\Gamma$, to be an arbitrary element of $\pi^{-1}(b)\cap A_\Gamma$.
Note that $\pi^{-1}(b)\cap A_\Gamma$ is non-empty (so that this choice is possible) as $b\in B_\Gamma\neq \varnothing$,
and, $B_\Gamma=\pi(A_\Gamma)$ by Lemma~\ref{common}.
The map $\alpha:B\rightarrow A$ thus obtained is clearly a~linear splitting of~$\pi$.
For any $i\in I$ consider any $b\in B_i$.
Then, by Lemma~\ref{partlem}, $b\in\Span(\bigcup_{\Gamma\in 2^I  \, |
\,
i\in\Gamma}\mathcal{B}_\Gamma)$, and hence
\begin{gather*}
\alpha(b)\in\sum\limits_{\Gamma\in 2^I
\,
|
\,
i\in\Gamma}\sum\limits_{b'\in \mathcal{B}_\Gamma}\big(\pi^{-1}(b')\cap A_\Gamma\big) \subseteq\sum\limits_{\Gamma\in
2^I
\,
|
\,
i\in\Gamma}A_\Gamma\subseteq A_i.\tag*{\qed}
\end{gather*}
\renewcommand{\qed}{}
\end{proof}

Finally, we argue that we can generate a~colinear splitting (with appropriate properties) from the linear one on the
coaction invariant subalgebra:

\begin{lemma}
\label{alphaexists}
Let~$A$ be a~principal $H$-comodule algebra, let $\pi:A\rightarrow B$ be an $H$-comodule algebra surjection, and let
$\{A_i\}_{i\in I}$ be a~finite family of ideals in~$A$ which are subcomodules, such that $\{A_i\}_{i\in
I}\cup\{\ker\pi\}$ generates a~distributive lattice.
Define for all $i\in I$
\begin{gather*}
A_i^{\co H}:=A_i\cap A^{\co H},
\qquad
B_i:=\pi(A_i),
\qquad
B_i^{\co H}:=B^{\co H}\cap B_i.
\end{gather*}
Suppose that there exists a~linear map $\alpha^{\co H}:B^{\co H}\rightarrow A^{\co H}$ such that
\begin{gather*}
\pi\circ\alpha^{\co H}=\mathrm{id}_{B^{\co H}},
\qquad
\alpha^{\co H}\big(B^{\co H}_i\big)\subseteq A_i^{\co H},
\qquad
\text{for all}
\quad
i\in I.
\end{gather*}
Let $\ell:H\rightarrow A\otimes A$ be a~strong connection on~$A$.
Then the following formula
\begin{gather*}
\alpha: \ B\longrightarrow A,
\qquad
b\longmapsto \alpha^{\co H}\big(b\sw{0}\pi\big(\ls{b\sw{1}}\big)\big)\rs{b\sw{1}}
\end{gather*}
defines a~right $H$-colinear map satisfying
\begin{gather*}
\pi\circ\alpha=\mathrm{id}_{B},
\qquad
\alpha(B_i)\subseteq A_i,
\qquad
\text{for all}
\quad
i\in I.
\end{gather*}
\end{lemma}
\begin{proof}
The fact that~$\alpha$ def\/ined above is a~colinear splitting of~$\pi$ follows immediately from the proof of~\cite[Lemma~3.1]{hkmz11}.
It remains to show that $\alpha(B_i)\subseteq A_i$ for all $i\in I$.
Indeed, let $b\in B_i$.
Because of the left colinearity of~$\ell$ (equation~\eqref{el-lcolinear}) it follows easily that
$b\sw{0}\pi(\ls{b\sw{1}})\otimes\rs{b\sw{1}}\in B^{\co H}\otimes A$ (cf.\
proof of~\cite[Lemma~3.1]{hkmz11}), and because $B_i$ is an ideal and also a~right $H$-subcomodule, it follows also that
$b\sw{0}\pi(\ls{b\sw{1}})\otimes\rs{b\sw{1}}\in B_i\otimes A$, hence $b\sw{0}\pi(\ls{b\sw{1}})\otimes\rs{b\sw{1}}\in B^{\co H}_i\otimes A$.
Therefore
\begin{gather*}
\alpha(b)=\alpha^{\co H}\big(b\sw{0}\pi\big(\ls{b\sw{1}}\big)\big)\rs{b\sw{1}} \in \alpha^{\co H}\big(B_i^{\co
H}\big)A\subseteq A^{\co H}_iA\subseteq A_i.\tag*{\qed}
\end{gather*}
\renewcommand{\qed}{}
\end{proof}

Let us now put together all the steps needed to construct a~strong connection on a~piecewise principal comodule algebra
using the results presented in this paper.
Let~$H$ be a~co-commutative Hopf algebra.
Suppose that~$P$ is an $H$-comodule algebra which is piecewise principal with respect to a~f\/inite family
$\{\pi_i:P\rightarrow P_i\}_{i\in J}$ of surjective $H$-comodule algebra morphisms.
Then by~\cite[Corollary~3.7]{hkmz11} the family $\{\pi_i:P\rightarrow P_i\}_{i\in J}$ is a~covering.
Hence, by~\cite[Proposition~3]{CM2000}, it follows that~$P$ is isomorphic with the multipullback $P^\pi$, where the
gluing morphisms are def\/ined, for all $i,j\in J$, $i\neq j$,~by
\begin{gather}
\label{CanFamily}
\pi^i_j: \ P_i\longrightarrow P_{ij}:=P/(\ker\pi_i+\ker\pi_j),
\qquad
\pi_i(p)\longmapsto p+\ker\pi_i+\ker\pi_j.
\end{gather}
Obviously $\pi^i_j$'s are surjective.
Note that $\ker\pi^i_j=\pi_i(\ker\pi_j)$ for all $i,j\in J$, $i\neq j$.
Hence, as $\ker\pi_i$'s generate a~distributive lattice of ideals (see above), it follows by Lemma~\ref{common} that,
for all $i\in J$, also $\ker\pi^i_j$'s generate a~distributive lattice of ideals.
It is also immediate (see, e.g.,~\cite[Remark~2]{CM2000}) that the family~\eqref{CanFamily} satisf\/ies the cocycle
condition (Def\/inition~\ref{cocycle-def}).
Before we can use Theorem~\ref{main-split} to construct splittings of $\pi_i$'s needed by Theorem~\ref{Main} we still
need to construct splittings $\alpha^i_j,\beta^i_j:P_{ij}\rightarrow P_i$, $i,j\in J$,
$j\neq i$ of $\pi^i_j$'s with appropriate properties.
The $H$-colinear splittings $\beta^i_j:P_{ij}\rightarrow P_i$ can be constructed using~\cite[Lemma~3.1]{hkmz11} as all
the $P_{i}$'s are principal by assumption.
Similarly, using Lemma~\ref{alphaexists} we can construct $H$-colinear splittings $\alpha^i_j:P_{ij}\rightarrow P_i$ of
$\pi^i_j$'s satisfying $\alpha^i_j(\pi^i_j(\ker{}\pi^i_k)\subseteq\ker\pi^i_k$ for all distinct $i,j,k\in J$ from
a~family of linear maps $\big(\alpha^{\co H}\big)^i_j:P_i^{\co H}\rightarrow P_{ij}^{\co H}$, $i,j\in J$, $i\neq j$ satisfying
\begin{gather}
\label{props2alph}
\big(\pi^{\co H}\big)^i_j\circ\big(\alpha^{\co H}\big)^i_j=\mathrm{id}_{P_{ij}^{\co H}},
\qquad
\big(\alpha^{\co H}\big)^i_j\big(\big(\pi^{\co H}\big)^i_j\big(\ker\big(\pi^{\co H}\big)^i_k\big)\big)\subseteq\ker\big(\pi^{\co H}\big)^i_k,
\end{gather}
where we denoted by $\big(\pi^{\co H}\big)^i_j:P_i^{\co H}\rightarrow P_{ij}^{\co H}$, $i,j\in J$, $i\neq j$
the appropriate restrictions of $\pi^i_j$'s to the coaction invariant subalgebras.
Note that by~\cite[Lemma~3.1]{hkmz11} all the $(\pi^{\co H})^i_j$'s are surjective.
Then Lemma~\ref{splitcoinv} gives a~semi-explicit construction of $(\alpha^{\co H})^i_j$'s satisfying properties~\eqref{props2alph}.

Using $\alpha^i_j$'s and $\beta^i_j$'s we can now construct $H$-colinear and unital splittings of $\pi_i$'s using
Theo\-rem~\ref{main-split} and utilize them in the explicit construction of a~strong connection given~by
Theorem~\ref{Main}.

\section{Example}\label{example-section}

In~\cite{reduct} a~new non-commutative real projective space $\mathbb{R} P^{2}_{\mathcal{T}}$ and a~non-commutative
sphere $S^2_{\mathbb{R}\mathcal{T}}$ were introduced, by def\/ining $\qcp 2$ and ${C(S^2_{\mathbb{R}\mathcal{T}})}$ as
a~particular triple pullbacks of, respectively, three copies of the Toeplitz algebra ${\mathcal{T}}$ and the tensor
product ${\mathcal{T}}\otimes C({\mathbb Z_2})$.
The algebra ${C(S^2_{\mathbb{R}\mathcal{T}})}$ has a~natural (component-wise) diagonal coaction of the Hopf algebra
$C({\mathbb Z_2})$, and it was proven in~\cite{reduct} that the subspace of invariants of this coaction is isomorphic
with $\qcp 2$.
Moreover, it was demonstrated that ${C(S^2_{\mathbb{R}\mathcal{T}})}$ is a~piecewise principal (hence principal)
$C({\mathbb Z_2})$-comodule algebra.
However, the paper~\cite{reduct} does not present an explicit formula for a~strong connection.
Because $C({\mathbb Z_2})$ is co-commutative and ${C(S^2_{\mathbb{R}\mathcal{T}})}$ is def\/ined as a~triple pullback
algebra, our main result is applicable.
In this section we will present the comparison of computations of a~strong connection on
${C(S^2_{\mathbb{R}\mathcal{T}})}$ using two methods: the f\/irst one uses the strong connection formula
from~\cite{hkmz11} and the other one uses Theorem~\ref{Main}.
The reader will see that, while application of the formula from~\cite{hkmz11} is trivial in case of double pullbacks,
already for triple pullbacks the computations becomes fairly unmanageable.
Also note that, in many cases, the values of strong connection formula on generators of the Hopf algebra are easily
guessable, and then the values on arbitrary Hopf algebra elements can be computed using well known recursive formula.
Here the Hopf algebra $C({\mathbb Z_2})$ has linear basis consisting of $1$ and~$u$, where~$u$ is the single generator
such that $u^2=1$, so that it suf\/f\/ices to f\/ind the value of a~strong connection on~$u$ without any need for recursion.
However, guessing the value of a~strong connection on~$u$ is nigh impossible.

We will start with recalling the def\/inition of the comodule algebra ${C(S^2_{\mathbb{R}\mathcal{T}})}$.
Our presentation will be very brief (mostly lifted from~\cite{reduct}), though suf\/f\/icient to understand what follows,
and will hardly include any geometric intuitions behind ${C(S^2_{\mathbb{R}\mathcal{T}})}$.
Also, because the def\/inition of~$\qcp 2$ is irrelevant for the strong connection computation, we omit it entirely.
Therefore, the reader is recommended to read the full account from~\cite{reduct}.

\subsection{A pullback quantum sphere}

We consider the Toeplitz algebra~$\mathcal{T}$ as the universal $C^*$-algebra generated by an isometry~$s$, and the
symbol map given by the assignment $\sigma\colon\mathcal{T}\ni s\mapsto \widetilde{u}\in C(S^1)$, where $\widetilde{u}$
is the unitary function generating $C(S^1)$.
The following two maps
\begin{gather*}
{\mathbb Z_2}\times I\ni(k,t)\stackrel{\delta_1}{\longmapsto} e^{i\pi\big(\frac{1}{4}kt+\frac{1}{2}k+\frac{3}{2}\big)}\in S^1,
\qquad
I\times{\mathbb Z_2}\ni(t,k)\stackrel{\delta_2}{\longmapsto} e^{i\pi\big({-}\frac{1}{4}kt-\frac{1}{2}k+1\big)}\in S^1,
\end{gather*}
and their pullbacks
\begin{gather*}
\delta^*_1: \  C\big(S^1\big)\longrightarrow C({\mathbb Z_2})\otimes C(I),
\qquad
\delta^*_2\colon C\big(S^1\big)\longrightarrow C(I)\otimes C({\mathbb Z_2}).
\end{gather*}
feature prominently in the def\/inition of ${C(S^2_{\mathbb{R}\mathcal{T}})}$.
We will denote for brevity $\sigma_i:=\delta^*_i\circ\sigma$, $i=1,2$.
The def\/initions of the $\delta_i$'s seem completely arbitrary.
In fact, as shown on the picture~\cite{reduct} below, each of these maps is meant as the parametrisation
of two appropriate quarters of $S^1$
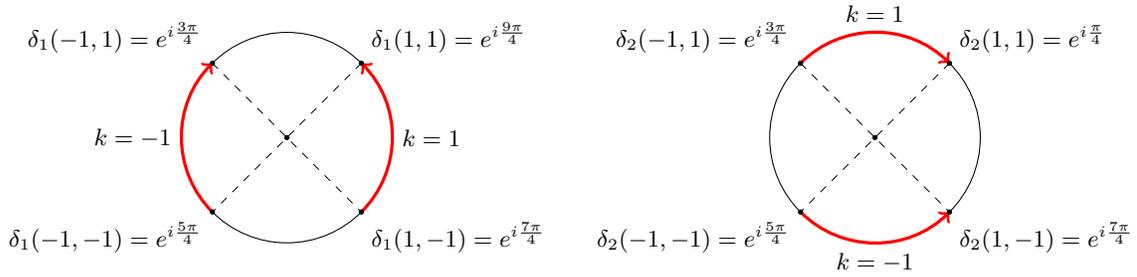
\begin{center}
\begin{tabular}{cc}
\begin{tikzpicture}[scale=0.7]
\fill (0,0) circle (0.5mm);
\draw (0,0) circle (2cm);
\draw[very thin, dashed] (45:2cm)--(225:2cm);
\draw[very thin, dashed] (135:2cm)--(315:2cm);
\draw[->,very thick,color=red] (225:2cm) arc (225:135:2cm);
\draw[->,very thick,color=red,rotate=-45] (0:2cm) arc (0:90:2cm);
\fill (45:2cm) circle (0.5mm) (135:2cm) circle (0.5mm)(225:2cm) circle (0.5mm) (315:2cm) circle (0.5mm);
\draw (-2cm,0pt) node[anchor=east]{{\footnotesize $k=-1$}};
\draw (2cm,0pt) node[anchor=west]{{\footnotesize $k=1$}};
\draw (0pt,2cm) node[anchor=south]{{\footnotesize $\phantom{k=1}$}};
\draw (0pt,-2cm) node[anchor=north]{{\footnotesize $\phantom{k=-1}$}};
\draw (45:2cm) node[anchor=south west]{{\footnotesize $\delta_1(1,1)=e^{i\frac{9\pi}{4}}$}};
\draw (315:2cm) node[anchor=north west]{{\footnotesize $\delta_1(1,-1)=e^{i\frac{7\pi}{4}}$}};
\draw (135:2cm) node[anchor=south east]{{\footnotesize $\delta_1(-1,1)=e^{i\frac{3\pi}{4}}$}};
\draw (225:2cm) node[anchor=north east]{{\footnotesize $\delta_1(-1,-1)=e^{i\frac{5\pi}{4}}$}};
\end{tikzpicture}
& \begin{tikzpicture}[scale=0.7]
\fill (0,0) circle (0.5mm);
\draw (0,0) circle (2cm);
\draw[very thin, dashed] (45:2cm)--(225:2cm);
\draw[very thin, dashed] (135:2cm)--(315:2cm);
\draw[->,very thick,color=red] (135:2cm) arc (135:45:2cm);
\draw[->,very thick,color=red] (225:2cm) arc (225:315:2cm);
\fill (45:2cm) circle (0.5mm) (135:2cm) circle (0.5mm)(225:2cm) circle (0.5mm) (315:2cm) circle (0.5mm);
\draw (0pt,2cm) node[anchor=south]{{\footnotesize $k=1$}};
\draw (0pt,-2cm) node[anchor=north]{{\footnotesize $k=-1$}};
\draw (45:2cm) node[anchor=south west]{{\footnotesize $\delta_2(1,1)=e^{i\frac{\pi}{4}}$}};
\draw (315:2cm) node[anchor=north west]{{\footnotesize $\delta_2(1,-1)=e^{i\frac{7\pi}{4}}$}};
\draw (135:2cm) node[anchor=south east]{{\footnotesize $\delta_2(-1,1)=e^{i\frac{3\pi}{4}}$}};
\draw (225:2cm) node[anchor=north east]{{\footnotesize $\delta_2(-1,-1)=e^{i\frac{5\pi}{4}}$}};
\end{tikzpicture} 
\end{tabular}
\end{center}

We view $S^1$ and~$I$ as ${\mathbb Z_2}$-spaces via multiplication by $\pm 1$.
Then ${\mathbb Z_2}\times I$ and $I\times {\mathbb Z_2}$ are ${\mathbb Z_2}$-spaces with the diagonal action.
Accordingly, $C(I)$, $C(S^1)$, $C({\mathbb Z_2})\otimes C(I)$ and $C(I)\otimes C({\mathbb Z_2})$ are right $C({\mathbb
Z_2})$-comodule algebras with coactions given by the pullbacks of respective ${\mathbb Z_2}$-actions.
Denote by~$u$ the generator $C({\mathbb Z_2})$ given by $u(\pm 1):=\pm 1$.
Then the assignment $s\mapsto s\otimes u$ makes ${\mathcal{T}}$ a~$C({\mathbb Z_2})$-comodule algebra.
(This coaction corresponds to the ${\mathbb Z_2}$-action given by $\alpha_{-1}^{\mathcal{T}}(s)=-s$.) It is easy to
verify that the maps $\delta_i$, $i=1,2$, are ${\mathbb Z_2}$-equivariant, so that their pullbacks $\delta^*_i$'s are
right $C({\mathbb Z_2})$-comodule maps.
Also, since the symbol map~$\sigma$ is a~right $C({\mathbb Z_2})$-comodule map, so are $\sigma_i$'s.

The construction of ${C(S^2_{\mathbb{R}\mathcal{T}})}$ can be seen as the quantum version of constructing the
topological 2-sphere by assembling three pairs of squares to the boundary of a~cube.
In the quantum version the algebra ${\mathcal{T}}\otimes C({\mathbb Z_2})$ replaces the pair of squares.
Explicitly, the algebra ${C(S^2_{\mathbb{R}\mathcal{T}})}$ is def\/ined in~\cite{reduct} to be the following triple
pullback of three copies of ${\mathcal{T}}\otimes C({\mathbb Z_2})$
\begin{gather*}
\xymatrix{
{\mathcal{T}}_0\otimes C({\mathbb Z_2})
\ar[d]_{\sigma_1\otimes\mathrm{id}} &
{\mathcal{T}}_1\otimes C({\mathbb Z_2})
\ar[d]^{\sigma_1\otimes\mathrm{id}}
\\
C({\mathbb Z_2})\otimes C(I)\otimes C({\mathbb Z_2}) &
\ar[l]^{\Phi_{01}} C({\mathbb Z_2})\otimes C(I)\otimes C({\mathbb Z_2}),
}
\\
\xymatrix{
{\mathcal{T}}_0\otimes C({\mathbb Z_2})
\ar[d]_{\sigma_2\otimes\mathrm{id}} &
{\mathcal{T}}_2\otimes C({\mathbb Z_2})
\ar[d]^{\sigma_1\otimes\mathrm{id}}
\\
C(I)\otimes C({\mathbb Z_2})\otimes C({\mathbb Z_2}) &
\ar[l]^{\Phi_{02}} C({\mathbb Z_2})\otimes C(I)\otimes C({\mathbb Z_2}),
}
\\
\xymatrix{
{\mathcal{T}}_1\otimes C({\mathbb Z_2})
\ar[d]_{\sigma_2\otimes\mathrm{id}} &
{\mathcal{T}}_2\otimes C({\mathbb Z_2})
\ar[d]^{\sigma_2\otimes\mathrm{id}}
\\
C(I)\otimes C({\mathbb Z_2})\otimes C({\mathbb Z_2}) &
\ar[l]^{\Phi_{12}} C(I)\otimes C({\mathbb Z_2})\otimes C({\mathbb Z_2}),
}
\end{gather*}
where the isomorphisms $\Phi_{ij}$ are def\/ined by the following formulas, for all $h,k\in C({\mathbb Z_2})$ and $p\in C(I)$
\begin{gather*}
\Phi_{01}(h\otimes p\otimes k) :=k\otimes p\otimes h,
\qquad
\Phi_{02}(h\otimes p\otimes k) :=p\otimes k\otimes h,
\\
\Phi_{12}(p\otimes h\otimes k) :=p\otimes k\otimes h.
\end{gather*}
We view the algebras ${\mathcal{T}}\otimes C({\mathbb Z_2})$, $C(I)\otimes C({\mathbb Z_2})\otimes C({\mathbb Z_2})$ and
$ C({\mathbb Z_2})\otimes C(I)\otimes C({\mathbb Z_2})$ as right $C({\mathbb Z_2})$-comodules with the diagonal
$C({\mathbb Z_2})$-coaction.
The coaction of $C({\mathbb Z_2})$ is def\/ined on ${C(S^2_{\mathbb{R}\mathcal{T}})}$ componentwise.

\subsection{Construction of certain auxilliary elements}

Both constructions of strong connections will require the existence of elements $\phi_1\in\sigma_1^{-1}(u\otimes
1_{C(I)})\subseteq{\mathcal{T}}$, $\phi_2\in\sigma_2^{-1}(1_{C(I)}\otimes u)\subseteq{\mathcal{T}}$ with certain
additional properties.
These elements will play the crucial role in the construction of appropriate splittings required by both methods.
More explicitly, we have the following:

\begin{lemma}
There exist elements $\phi_1,\phi_2\in {\mathcal{T}}$ satisfying
\begin{subequations}
\label{phiprop2}
\begin{alignat}{3}
& \rho(\phi_1)=\phi_1\otimes u,
\qquad &&
\rho(\phi_2)=\phi_2\otimes u, &
\label{phiprop2a}
\\
& \sigma_1(\phi_1)=u\otimes 1_{C(I)},
\qquad &&
\sigma_2(\phi_1)=\imath_I\otimes 1_{C({\mathbb Z_2})},&
\label{phiprop2ba}
\\
& \sigma_2(\phi_2)=1_{C(I)}\otimes u,
\qquad &&
\sigma_1(\phi_2)=1_{C({\mathbb Z_2})}\otimes \imath_I,&
\label{phiprop2b}
\\
& \big(1-\phi_2^2\big)\big(1-\phi_1^2\big)\neq0,\qquad &&&
\label{phiprop2c}
\end{alignat}
\end{subequations}
where $\imath_I\in C(I)$ is an an identity map $\imath_I(t)=t$ and $\rho:{\mathcal{T}}\rightarrow {\mathcal{T}}\otimes
C({\mathbb Z_2})$ is a~right coaction.
\end{lemma}
\begin{proof}

First we def\/ine auxiliary maps $\hat\phi_1,\hat\phi_2\in C(S^1)$ by the formulae
\begin{gather}
\hat\phi_1\big(e^{i\theta}\big):=
\begin{cases}
2-\frac{4}{\pi}\theta & \text{if} \quad \theta\in[\frac{\pi}{4},\frac{3\pi}{4}],
\\
-1 & \text{if} \quad  \theta\in[\frac{3\pi}{4},\frac{5\pi}{4}],
\\
\frac{4}{\pi}\theta-6 &\text{if} \quad  \theta\in[\frac{5\pi}{4},\frac{7\pi}{4}],
\\
1 & \text{if} \quad  \theta\in[\frac{7\pi}{4},\frac{9\pi}{4}],
\end{cases}
\qquad
\hat\phi_2\big(e^{i\theta}\big):=
\begin{cases}
1 & \text{if} \quad \theta\in[\frac{\pi}{4},\frac{3\pi}{4}],
\\
4-\frac{4}{\pi}\theta & \text{if} \quad  \theta\in[\frac{3\pi}{4},\frac{5\pi}{4}],
\\
-1 &\text{if} \quad  \theta\in[\frac{5\pi}{4},\frac{7\pi}{4}],
\\
\frac{4}{\pi}\theta-8 &\text{if} \quad  \theta\in[\frac{7\pi}{4},\frac{9\pi}{4}].
\end{cases}
\label{phidefeq}
\end{gather}
One immediately verif\/ies that
\begin{gather}
\label{line2}
\hat\phi_1,\hat\phi_2: \ S^1\longrightarrow [-1,1],
\qquad
\hat\phi_1(-z)=-\hat\phi_1(z),
\qquad
\hat\phi_2(-z)=-\hat\phi_2(z).
\end{gather}
(i.e., $\rho(\hat\phi_1)=\hat\phi_1\otimes u$ and $\rho(\hat\phi_2)=\hat\phi_2\otimes u$) and that
\begin{alignat}{3}
& \hat\phi_1\circ\delta_1=u\otimes 1_{C(I)},
\qquad&&
\hat\phi_2\circ\delta_2=1_{C(I)}\otimes u,&
\nonumber
\\
&\hat\phi_1\circ\delta_2=\imath_I\otimes 1_{C({\mathbb Z_2})},
\qquad&&
\hat\phi_2\circ\delta_1=1_{C({\mathbb Z_2})}\otimes \imath_I.&
\label{phiprop}
\end{alignat}
Using equation~\eqref{line2} and the standard properties of comodules, one proves that because the symbol map~$\sigma$
is a~surjective right $C({\mathbb Z_2})$ comodule map and~$u$ is grouplike, we can choose elements
$\phi_i\in{\mathcal{T}}$, $i=1,2$, such that $\sigma(\phi_i)=\hat\phi_i$, and $\rho(\phi_i)=\phi_i\otimes u$, thus
verifying the properties~\eqref{phiprop2a}.
That thus chosen elements $\phi_1$, $\phi_2$ satisfy the properties~\eqref{phiprop2ba} and~\eqref{phiprop2b} follows
immediately from equations~\eqref{phiprop}.

The last condition of the lemma is an easy consequence of the properties of the representation of a~Toeplitz algebra on
a~Bergman space (see, e.g.,~\cite[Theorem~2.8.2]{vasil}).
However, we provide an alternative elementary proof to make the presentation self-contained.
Unfortunately, $\sigma((1-\phi_2^2)(1-\phi_1^2))=(1-\hat\phi_2^2)(1-\hat\phi_1^2)=0$, so we cannot prove that vector
$(1-\phi_2^2)(1-\phi_1^2)\in{\mathcal{T}}$ is nonzero by considering the properties of its image in $C(S^1)$
under~$\sigma$ and we must work directly in ${\mathcal{T}}$.
We will use the f\/lexibility af\/forded by the fact that conditions~\eqref{phiprop2a},~\eqref{phiprop2ba}
and~\eqref{phiprop2b} do not f\/ix completely elements $\phi_i\in{\mathcal{T}}$.
We will show that even if $(1-\phi_2^2)(1-\phi_1^2)=0$ for our initial choice of $\phi_i$'s, there exists a~family
$\{\phi_{2;t,n}\}_{t\in{\mathbb R},n\in{\mathbb N}}$ of deformations of $\phi_2$ such that the
conditions~\eqref{phiprop2a},~\eqref{phiprop2ba} and~\eqref{phiprop2b} are still satisf\/ied for all pairs
$(\phi_1,\phi_{2;t,n})$ and there exist $n\in{\mathbb N}$ and $t\in{\mathbb R}$ such that $(1-\phi_{2;t,n}^2)(1-\phi_1^2)\neq
0$.

Let~$z$ be a~partial isometry generating ${\mathcal{T}}$, and let $\rho:{\mathcal{T}}\rightarrow{\mathcal{T}}\otimes
C({\mathbb Z_2})$ be a~right $C({\mathbb Z_2})$-coaction.
Def\/ine, for all $n\in{\mathbb N}$, $t\in{\mathbb R}$
\begin{gather}
\phi_{2;t,n}:=\phi_2+tE_n,
\qquad
\text{where}
\qquad
E_n=z\left(z^n(z^*)^n-z^{n+2}(z^*)^{n+2}\right).
\label{phindef}
\end{gather}
Because $\rho(z)=z\otimes u$, we have $\rho(\phi_{2;t,n})=\phi_{2;t,n}\otimes u$ and because $\sigma(E_n)=0$, we have
$\sigma(\phi_{2;t,n})=\hat\phi_2$, hence all of the conditions~\eqref{phiprop2} are satisf\/ied, and for all~$t$ and~$n$
we can use $\phi_{2;t,n}$ instead of $\phi_2$ in the formula~\eqref{ellform3} def\/ining a~strong connection on
${C(S^2_{\mathbb{R}\mathcal{T}})}$.
Assume that $(1-\phi_{2;t,n}^2)(1-\phi_1^2)=0$ for all $t\in{\mathbb R}$ and $n\in{\mathbb N}$.
We will show that this assumption leads to contradiction.
Using equation~\eqref{phindef} elements $(1-\phi_{2;t,n}^2)(1-\phi_1^2)$ can be explicitly written as
\begin{gather*}
\big(1-\phi_2^2\big)\big(1-\phi_1^2\big)-(E_n\phi_2+\phi_2E_n)t-(E_n)^2\big(1-\phi_1^2\big)t^2.
\end{gather*}
If $(1-\phi_{2;t,n}^2)(1-\phi_1^2)=0$ for all $t\in{\mathbb R}$ and $n\in{\mathbb N}$ then the above polynomials in~$t$ are
identically zero for all $n\in{\mathbb N}$, which implies in particular that coef\/f\/icients at $t^2$ must be zero, i.e., that
\begin{gather}
\label{absurdcond}
E_n^2\big(1-\phi_1^2\big)=0,
\qquad
\text{for all}
\quad
n\in N.
\end{gather}

Consider now the faithful representation $R:{\mathcal{T}}\rightarrow\mathcal{H}$ of the Toeplitz algebra ${\mathcal{T}}$
on a~Hilbert space $\mathcal{H}$ spanned by an orthonormal basis $|n\rangle$, $n\in{\mathbb N}$, where the partial
isometry~$z$ is represented as a~right shift, i.e., $R(z)|n\rangle=|n+1\rangle$ for all $n\in{\mathbb N}$.
One easily proves that
\begin{gather}
\label{Erep}
R(E_n^2)|m\rangle=\delta_{m,n}|n+2\rangle,
\qquad
\text{for all}\quad  m,n\in{\mathbb N}.
\end{gather}
Equation~\eqref{absurdcond} implies that $R(E_n^2)R(1-\phi_1^2)\Psi=0$, for all $n\in{\mathbb N}$ and $\Psi\in\mathcal{H}$.
But then it follows from equation~\eqref{Erep} that $R(1-\phi_1^2)\Psi=0$ for all $\Psi\in\mathcal{H}$, i.e., that
$R(1-\phi_1^2)=0$.
But~$R$ is faithful, hence $(1-\phi_1^2)=0$.
On the other hand, $\sigma(1-\phi_1^2)=1-\hat\phi_1^2\neq 0$.
Hence we reached contradiction.
It follows that we can choose $\phi_1$ and $\phi_2$ so that all conditions~\eqref{phiprop2} are satisf\/ied.
\end{proof}

\subsection{A~strong connection. Method~I}

In this subsection we construct a~strong connection on the $C({\mathbb Z_2})$-comodule
algebra~${C(S^2_{\mathbb{R}\mathcal{T}})}$ by repeated application of the formula stated in the proof of~\cite[Lemma~3.2]{hkmz11}.
Let~$P$ be a~f\/ibre product of $P_1\stackrel{\pi^1_2}{\longrightarrow}P_{12}\stackrel{\pi^1_2}{\longleftarrow}P_2$ in the
category of right $H$-comodule algebras.
Assume that the maps $\pi^i_j$ are surjective and that $\ell_i:H\rightarrow P_i\otimes P_i$, $i=1,2$, are strong connections.
Then the formula \cite[Lemma~3.2]{hkmz11}
\begin{gather}
\ell(h) =\big(\ell_1(h){}^{\langle 1\rangle},f^1_2(\ell_1(h){}^{\langle
1\rangle})\big)\otimes\big(\ell_1(h){}^{\langle 2\rangle},f^1_2\big(\ell_1(h){}^{\langle 2\rangle}\big)\big)
\nonumber
\\
\phantom{\ell(h) =}
{}+\big(0,\big(\varepsilon(h\sw{1})-f^1_2(\ell_1(h\sw{1}){}^{\langle 1\rangle})f^1_2(\ell_1(h\sw{1}){}^{\langle
2\rangle})\big)
\ell_2(h\sw{2}){}^{\langle 1\rangle} \big)
\nonumber
\\
\phantom{\ell(h) =+}{}
\otimes\big(f^2_1\big(\ell_2(h\sw{2}){}^{\langle 2\rangle}\big),\ell_2(h\sw{2}){}^{\langle 2\rangle}\big)
\label{stronglue}
\end{gather}
def\/ines a~strong connection $\ell:H\rightarrow P\otimes P$.
Here $f^i_j:=\mu_j^i\circ\pi^i_j$ and $\mu_j^i$ is any unital colinear splitting of $\pi^j_i$, $i\neq j$.
Note also that we use the convention that, if $x{}^{\langle 1\rangle}\otimes x{}^{\langle
2\rangle}:=\sum\limits_ix_i\otimes y_i$, then $(x{}^{\langle 1\rangle},x{}^{\langle 1\rangle})\otimes (x{}^{\langle
2\rangle},x{}^{\langle 2\rangle}):=\sum\limits_i(x_i,x_i)\otimes (y_i,y_i)$, and similarly for coproducts.
Observe that for $C({\mathbb Z_2})$-comodule algebras it is enough to compute the value of a~strong connection for
$h=u$, where~$u$ is the group-like generator of $C({\mathbb Z_2})$ because strong connections are unital and linear,
i.e., it is suf\/f\/icient to use the following equation
\begin{gather}
\ell(u) =\big(\ell_1(u){}^{\langle 1\rangle},f^1_2\big(\ell_1(u){}^{\langle
1\rangle}\big)\big)\otimes\big(\ell_1(u){}^{\langle 2\rangle},f^1_2\big(\ell_1(u){}^{\langle 2\rangle}\big)\big)
\nonumber
\\
\phantom{\ell(u) =}
{}+\big(0,\big(1-f^1_2\big(\ell_1(u){}^{\langle 1\rangle}\big)f^1_2\big(\ell_1(u){}^{\langle 2\rangle}\big)\big)
\ell_2(u){}^{\langle 1\rangle} \big)\otimes\big(f^2_1\big(\ell_2(u){}^{\langle 2\rangle}\big),\ell_2(u){}^{\langle 2\rangle}\big).
\label{stronglue2}
\end{gather}
Note that it is suf\/f\/icient to know the values $f^i_j(x)$ only for a~set of elements $x\in P_j$ which actually appear in
the above formula and which (because of bi-colinearity of strong connections) can be assumed to be linearly independent
and satisfy $\rho(x)=x\otimes u$, i.e., one needs only to solve the following equations with unknowns $f^i_j(x)\in P_j$
(where~$\rho$ denotes the coaction)
\begin{gather}
\label{fijeqsys}
\rho\big(f^i_j(x)\big)=f^i_j(x)\otimes u,
\qquad
\pi^j_i\big(f^i_j(x)\big)=\pi^i_j(x).
\end{gather}

As the formula~\eqref{stronglue} assumes the comodule algebra to be presented as the ordinary (double) pullback, we need
to convert the triple-pullback def\/ining ${C(S^2_{\mathbb{R}\mathcal{T}})}$ to an iterated pullback and apply the formula
recursively.
Since all the maps $C(\proj R 2)\to \mathcal{T}_i\otimes C({\mathbb Z_2})$ are surjective~\cite{reduct}, we can
apply~\cite[Lemma~0.2 and Proposition~1.3]{r-j} to present ${C(S^2_{\mathbb{R}\mathcal{T}})}$ as a~desired iterated pullback
\begin{gather}
\label{iterpullback}
\begin{split}
\xymatrix@C=65pt@R=40pt@!0{
&&& {C(S^2_{\mathbb{R}\mathcal{T}})}\ar[dll]\ar[drr]&&
\\
& P_1\ar[dr]\ar[drrr]^{\beta_2}\ar[dl] &&&& (\mathcal T\otimes C({\mathbb Z_2}))_2\ar[dl]_{\beta_1}
\\
(\mathcal T\otimes C({\mathbb Z_2}))_0\ar@{.>}[drrr]
\ar[dr]_{\alpha_1}&&(\mathcal T\otimes C({\mathbb Z_2}))_1\ar[dl]^{\alpha_2}\ar@{.>}[drrr]&&\lim P_{12}\ar[dl]\ar[dr]&
\\
& C({\mathbb Z_2})\otimes\mathcal T\otimes C({\mathbb Z_2}) && \mathcal T\otimes C({\mathbb Z_2})\otimes C({\mathbb
Z_2})\ar[dr] && \mathcal T\otimes C({\mathbb Z_2})\otimes C({\mathbb Z_2}).\ar[dl]
\\
&&&& C({\mathbb Z_2})\otimes C({\mathbb Z_2})\otimes C({\mathbb Z_2})&}
\end{split}\!\!\!\!\!
\end{gather}
Here
\begin{alignat*}{3}
& \alpha_1 =\sigma_1\otimes \mathrm{id},
\qquad&&
\beta_1(x) =((\Phi_{02}\circ(\sigma_1\otimes\mathrm{id}))(x),(\Phi_{12}\circ(\sigma_2\otimes\mathrm{id}))(x)),&
\\
& \alpha_2 =\Phi_{01}\circ(\sigma_1\otimes\mathrm{id}),
\qquad &&
\beta_2(x,y) =((\sigma_2\otimes\mathrm{id})(x),(\sigma_2\otimes\mathrm{id})(y)).&
\end{alignat*}

We will f\/irst compute a~strong connection $\ell_{01}:C({\mathbb Z_2})\rightarrow P_1\otimes P_1$ on $P_1$~-- the f\/iber
product of ${\mathcal{T}}_0\otimes C({\mathbb Z_2})$ and ${\mathcal{T}}_1\otimes C({\mathbb Z_2})$
(see~\eqref{iterpullback}).
We use the particular choice of the strong connections $\ell_0$ and $\ell_1$ on trivial pieces
${\mathcal{T}}_0\otimes C({\mathbb Z_2})$ and ${\mathcal{T}}_1\otimes C({\mathbb Z_2})$ given~by
\begin{gather*}
\ell_0(u)=(1\otimes u)\otimes (1\otimes u),
\qquad
\ell_1(u)=(1\otimes u)\otimes (1\otimes u).
\end{gather*}
Substituting the above formulae in~\eqref{stronglue2} yields a~glued strong connection on $P_1$
\begin{gather}
\ell_{01}(u) =\big(1\otimes u,f^0_1(1\otimes u)\big)\otimes\big(1\otimes u,f^0_1(1\otimes u)\big)
\nonumber
\\
\hphantom{\ell_{01}(u) =}
{}+\big(0,\big(1\otimes 1-f^0_1(1\otimes u)f^0_1(1\otimes u)\big)(1\otimes u) \big)\otimes\big(f^1_0(1\otimes u),1\otimes u
\big).
\label{strong01}
\end{gather}
Let us write for brevity $a:=f^0_1(1\otimes u)$, $b:=f^1_0(1\otimes u)$.
By diagram~\eqref{iterpullback} and equation~\eqref{fijeqsys} elements $a,b\in{\mathcal{T}}\otimes C({\mathbb Z_2})$ are any
solutions to the following equations
\begin{gather*}
\rho(a) =a\otimes u,  (\sigma_1\otimes\mathrm{id})(1\otimes u) =(\Phi_{01}\circ(\sigma_1\otimes\mathrm{id}))(a),
\nonumber
\\
\rho(b) =b\otimes u,  (\sigma_1\otimes\mathrm{id})(b) =(\Phi_{01}\circ(\sigma_1\otimes\mathrm{id}))(1\otimes u).
\end{gather*}
Substituting the def\/inition of $\Phi_{01}$ simplif\/ies the above system of equations to
\begin{gather*}
\rho(a)=a\otimes u,
\qquad\!
\rho(b)=b\otimes u,
\qquad\!
(\sigma_1\otimes\mathrm{id})(a)=u\otimes 1\otimes1,
\qquad\!
(\sigma_1\otimes\mathrm{id})(b)=u\otimes 1\otimes1,
\end{gather*}
and it is easy to see that one of the solutions is
\begin{gather*}
a=\phi_1\otimes1,
\qquad
b=\phi_1\otimes1.
\end{gather*}
Here and in what follows $\phi_1$ and $\phi_2$ are elements of ${\mathcal{T}}$ satisfying all the conditions~\eqref{phiprop2}.
Substituting the above solution into~\eqref{strong01} yields the following strong connection on $P_1$
\begin{gather}
\label{strong01f}
\ell_{01}(u)=(1\otimes u,\phi_1\otimes 1)\otimes(1\otimes u,\phi_1\otimes 1) +(0,(1-\phi_1^2)\otimes u)\otimes
(\phi_1\otimes1,1\otimes u)
\end{gather}

Now, we apply the formula~\eqref{stronglue} to the second iterated pullback in the diagram~\eqref{iterpullback}
\begin{gather*}
\vcenter{\xymatrix@C=10pt{{\mathcal{T}}_2\otimes C({\mathbb Z_2})\ar[dr]_{\beta_1}&&P_1\ar[dl]^{\beta_2}
\\
&Q&}}
\qquad
\begin{array}{l}
\beta_1(x)=((\Phi_{02}\circ(\sigma_1\otimes\mathrm{id}))(x),(\Phi_{12}\circ(\sigma_2\otimes\mathrm{id}))(x)),
\\
\beta_2(x,y)=((\sigma_2\otimes\mathrm{id})(x),(\sigma_2\otimes\mathrm{id})(y)).
\end{array}
\end{gather*}
where $Q\subseteq (C(I)\otimes C({\mathbb Z_2})\otimes C({\mathbb Z_2}))\oplus (C(I)\otimes C({\mathbb Z_2})\otimes C({\mathbb Z_2}))$.
We choose the strong connection on $P_1$ given by equation~\eqref{strong01f}, and on ${\mathcal{T}}_2\otimes C({\mathbb Z_2})$ given~by
\begin{gather*}
\ell_2(u)=(1\otimes u)\otimes (1\otimes u).
\end{gather*}
Substituting these into formula~\eqref{stronglue2} yields a~strong connection on ${C(S^2_{\mathbb{R}\mathcal{T}})}$
\begin{gather}
\ell(u) =\big(\ell_2(u){}^{\langle 1\rangle},f^2_{01}\big(\ell_2(u){}^{\langle 1\rangle}\big)\big)\otimes\big(\ell_2(u){}^{\langle
2\rangle},f^2_{01}\big(\ell_2(u){}^{\langle 2\rangle}\big)\big)
\nonumber
\\
\hphantom{\ell(u)=}
{}+\big(0,\big((1\otimes1,1\otimes 1)-f^2_{01}(\ell_2(u){}^{\langle 1\rangle})f^2_{01}(\ell_2(u){}^{\langle
2\rangle})\big)\ell_{01}(u){}^{\langle 1\rangle} \big)
\nonumber
\\
\hphantom{\ell(u)=+}{}
\otimes\big(f^{01}_2\big(\ell_{01}(u){}^{\langle 2\rangle}\big),\ell_{01}(u){}^{\langle 2\rangle} \big)
\nonumber
\\
\hphantom{\ell(u)}
 =\big(1\otimes u,f^2_{01}(1\otimes u)\big)\otimes\big(1\otimes u,f^2_{01}(1\otimes u)\big)
\nonumber
\\
\hphantom{\ell(u)=}
{}+\big(0,\big((1\otimes1,1\otimes 1)-\big(f^2_{01}(1\otimes u)\big)^2\big)(1\otimes u,\phi_1\otimes 1) \big)
\nonumber
\\
\phantom{\ell(u)=+}
\otimes\big(f^{01}_2((1\otimes u,\phi_1\otimes 1)),(1\otimes u,\phi_1\otimes 1) \big)
\nonumber
\\
\phantom{\ell(u)=}
{}+\big(0,\big((1\otimes1,1\otimes 1)-\big(f^2_{01}(1\otimes u)\big)^2\big)\big(0,\big(1-\phi_1^2\big)\otimes u\big) \big)
\nonumber
\\
\hphantom{\ell(u)=+}{}
\otimes\big(f^{01}_2((\phi_1\otimes1,1\otimes u)),(\phi_1\otimes1,1\otimes u) \big),
\label{ellform1}
\end{gather}
where $f^2_{01}:{\mathcal{T}}\otimes C({\mathbb Z_2})\rightarrow P_1$, $f^{01}_2:P_1\rightarrow {\mathcal{T}}\otimes
C({\mathbb Z_2})$ are any linear, unital, right $C({\mathbb Z_2})$-comodule maps satisfying $\beta_2\circ
f^2_{01}=\beta_1$, $\beta_1\circ f^{01}_2=\beta_2$. Denote for brevity $(a_0,a_1):=f^2_{01}(1\otimes u)$,
$b:=f^{01}_2((1\otimes u,\phi_1\otimes 1))$, $c:=f^{01}_2((\phi_1\otimes1,1\otimes u))$.
It follows that we need to solve the following system of equations for $a_0$, $a_1$,~$b$,~$c$
\begin{alignat*}{3}
&\rho(a_0) =a_0\otimes u, \qquad&&  \big(\Phi_{02}\circ(\sigma_1\otimes\mathrm{id})\big)(1\otimes u) =(\sigma_2\otimes\mathrm{id})(a_0),&
\\
&\rho(a_1) =a_1\otimes u,  \qquad&&  \big(\Phi_{12}\circ(\sigma_2\otimes\mathrm{id})\big)(1\otimes u) =(\sigma_2\otimes\mathrm{id})(a_1),&
\\
&\rho(b) =b\otimes u,  \qquad&&  \big(\Phi_{02}\circ(\sigma_1\otimes\mathrm{id})\big)(b) =\sigma_2(1)\otimes u,&
\\
&&&
\big(\Phi_{12}\circ(\sigma_2\otimes\mathrm{id})\big)(b) =\sigma_2(\phi_1)\otimes 1_{C({\mathbb Z_2})},&
\\
& \rho(c) =c\otimes u,  \qquad&&  \big(\Phi_{02}\circ(\sigma_1\otimes\mathrm{id})\big)(c) =\sigma_2(\phi_1)\otimes 1_{C({\mathbb Z_2})},&
\\
&&&
\big(\Phi_{12}\circ(\sigma_2\otimes\mathrm{id})\big)(c) =\sigma_2(1)\otimes u.&
\end{alignat*}
Simplif\/ication of the right column of the above equations using equation~\eqref{phiprop2} yields
\begin{alignat*}{3}
&(\sigma_2\otimes\mathrm{id})(a_0) =1_{C(I)}\otimes u\otimes 1_{C({\mathbb Z_2})},
\qquad&&
(\sigma_2\otimes\mathrm{id})(a_1) =1_{C(I)}\otimes u\otimes 1_{C({\mathbb Z_2})},&
\\
&(\sigma_1\otimes\mathrm{id})(b) =u\otimes 1_{C(I)}\otimes 1_{C({\mathbb Z_2})},
\qquad&&
(\sigma_2\otimes\mathrm{id})(b) =\imath_I\otimes 1_{C({\mathbb Z_2})} \otimes 1_{C({\mathbb Z_2})},&
\\
&(\sigma_1\otimes\mathrm{id})(c) =1_{C({\mathbb Z_2})}\otimes \imath_I\otimes 1_{C({\mathbb Z_2})},
\qquad&&
(\sigma_2\otimes\mathrm{id})(c) =1_{C(I)}\otimes u\otimes 1_{C({\mathbb Z_2})}.&
\end{alignat*}
Using equation~\eqref{phiprop} again, one easily verif\/ies that one of the solutions can be given as
\begin{gather*}
(a_0,a_1)=\big(\phi_2\otimes 1_{C({\mathbb Z_2})},\phi_2\otimes 1_{C({\mathbb Z_2})}\big),
\qquad
b=\phi_1\otimes 1_{C({\mathbb Z_2})},
\qquad
c=\phi_2\otimes 1_{C({\mathbb Z_2})}.
\end{gather*}
Substituting this particular solution to the formula~\eqref{ellform1} for $\ell(u)$ yields
\begin{gather*}
\ell(u) =(1\otimes u,(\phi_2\otimes1,\phi_2\otimes 1))\otimes(1\otimes u,(\phi_2\otimes1,\phi_2\otimes 1))
\nonumber
\\
\hphantom{\ell(u) =}
{}+\big(0,\big((1\otimes1,1\otimes 1)-((\phi_2\otimes1,\phi_2\otimes 1))^2\big)(1\otimes u,\phi_1\otimes 1) \big)
\nonumber
\\
\hphantom{\ell(u) =+}{}
\otimes (\phi_1\otimes1,(1\otimes u,\phi_1\otimes 1) )
\nonumber
\\
\hphantom{\ell(u) =}
{}+\big(0, \big((1\otimes1,1\otimes 1)-((\phi_2\otimes1,\phi_2\otimes 1))^2\big)\big(0,\big(1-\phi_1^2\big)\otimes u\big) \big)
\nonumber
\\
\hphantom{\ell(u) =+}{}
\otimes (\phi_2\otimes1,(\phi_1\otimes1,1\otimes u)  ).
\end{gather*}
Simplifying, removing unnecessary parentheses and rearranging terms so that
${C(S^2_{\mathbb{R}\mathcal{T}})}\subseteq({\mathcal{T}}_0\otimes C({\mathbb Z_2}))\oplus ({\mathcal{T}}_1\otimes
C({\mathbb Z_2}))\oplus ({\mathcal{T}}_2\otimes C({\mathbb Z_2}))$ yields f\/inally
\begin{gather}
\ell(u) =(\phi_2\otimes1,\phi_2\otimes1,1\otimes u)\otimes(\phi_2\otimes1,\phi_2\otimes1,1\otimes u)
\nonumber
\\
\hphantom{\ell(u) =}
{}+ \big(\big(1-\phi_2^2\big)\otimes u,\big(1-\phi_2^2\big)\phi_1\otimes1,0\big) \otimes(1\otimes u,\phi_1\otimes1,\phi_1\otimes 1)
\nonumber
\\
\hphantom{\ell(u) =}
{}+ \big(0,\big(1-\phi_2^2\big)\big(1-\phi_1^2\big)\otimes u,0\big) \otimes(\phi_1\otimes1,1\otimes u,\phi_2\otimes 1).
\label{ellform3}
\end{gather}

Write $\ell(u)=\sum\limits_{i=1}^3l_i\otimes r_i$ where
\begin{alignat*}{3}
&l_1 =(\phi_2\otimes1,\phi_2\otimes1,1\otimes u),
\qquad&&
r_1 =(\phi_2\otimes1,\phi_2\otimes1,1\otimes u),&
\\
&l_2 =\big(\big(1-\phi_2^2\big)\otimes u,\big(1-\phi_2^2\big)\phi_1\otimes1,0\big),
\qquad&&
r_2 =(1\otimes u,\phi_1\otimes1,\phi_1\otimes 1),&
\\
& l_3 =\big(0,\big(1-\phi_2^2\big)\big(1-\phi_1^2\big)\otimes u,0\big),
\qquad&&
r_3 = (\phi_1\otimes1,1\otimes u,\phi_2\otimes 1).&
\end{alignat*}
According to~\cite[Theorem~3.1]{bh04} if both $\{l_1,l_2,l_3\}$ and $\{r_1,r_2,r_3\}$ are (separately) sets of linearly
independent vectors then
\begin{gather*}
p_{ij}:=r_il_j,
\qquad
p:=(p_{ij})\in M_3\big(C\big(\proj R2\big)\big),
\qquad
p^2=p,
\end{gather*}
is a~projector for an associated line bundle.
Hence, in order to use this result, we need to prove that zeros are the only solutions to equations
\begin{gather}
\sum\limits_{i=1}^3\alpha_ir_i=0,
\qquad
\sum\limits_{i=1}^3\beta_il_i=0.
\label{linindep1}
\end{gather}
The f\/irst \looseness=-1 of the above equalities implies that $\alpha_1\phi_2\otimes 1+\alpha_21\otimes u+\alpha_3\phi_1\otimes 1=0$,
hence immediately $\alpha_2=0$.
Because $\{\hat\phi_1,\hat\phi_2\}=\sigma(\{\phi_1,\phi_2\})$ (equation~\eqref{phiprop2}) are linearly independent in
$C(S^1)$, which can be checked easily by direct computation using equation~\eqref{phidefeq}, also $\{\phi_1,\phi_2\}$ must be
linearly independent in ${\mathcal{T}}$, hence $\alpha_1=\alpha_3=0$.
The second equality in equation~\eqref{linindep1} can be expanded as
\begin{gather}
\beta_1\phi_2\otimes 1+\beta_2\big(1-\phi^2_2\big)\otimes u =0,
\label{linindep0}
\\
\beta_1\phi_2\otimes 1+\beta_2\big(1-\phi^2_2\big)\phi_1\otimes 1+\beta_3\big(1-\phi_2^2\big)\big(1-\phi_1^2\big)\otimes u =0,
\label{linindep2}
\\
\beta_1 1\otimes u =0.
\label{linindep3}
\end{gather}
It follows immediately from equation~\eqref{linindep3} that $\beta_1=0$.
Then because $\sigma(1-\phi_2^2)=1-\hat\phi_2^2\neq 0$ (see equations~\eqref{phiprop2} and \eqref{phidefeq}) we have also
$1-\phi_2^2\neq 0$, and so, by equation~\eqref{linindep0}, $\beta_2=0$.
Finally, equation~\eqref{linindep2} and equation~\eqref{phiprop2c} implies that $\beta_3=0$.

\subsection{A~strong connection. Method~II}

In this subsection we construct a~strong connection on the $C({\mathbb Z_2})$-comodule
algebra~${C(S^2_{\mathbb{R}\mathcal{T}})}$ using the formula given in Theorem~\ref{Main}.
As in the previous subsection $\phi_1,\phi_2\in{\mathcal{T}}$ denote some chosen elements satisfying all the conditions in~\eqref{phiprop2}.
Also as before, the strong connections on the three copies of $C({\mathbb Z_2})$-comodule algebra (with diagonal coaction)
${\mathcal{T}}\otimes C({\mathbb Z_2})$ which are components of ${C(S^2_{\mathbb{R}\mathcal{T}})}$ are chosen as given by the formulas
\begin{gather*}
\ell_1(u)=\ell_2(u)=\ell_3(u)=1_{\mathcal{T}}\otimes u,
\qquad
\ell_1(1_{C({\mathbb Z_2})})=\ell_2(1_{C({\mathbb Z_2})})=\ell_3(1_{C({\mathbb Z_2})})=1_{\mathcal{T}}\otimes 1_{C({\mathbb Z_2})}.
\end{gather*}
In order to use the formula from Theorem~\ref{Main} we need the appropriate colinear and unital
splittings from the linear subspaces generated by the legs of $\ell_i$'s into ${C(S^2_{\mathbb{R}\mathcal{T}})}$.
The reader will easily verify that the maps $\alpha_i:\Span\{1_{\mathcal{T}}\otimes u,1_{\mathcal{T}}\otimes
1_{C({\mathbb Z_2})}\}\rightarrow{C(S^2_{\mathbb{R}\mathcal{T}})}$, $i=0,1,2$ def\/ined by setting
$\alpha_i(1_{\mathcal{T}}\otimes 1_{C({\mathbb Z_2})}):=1_{{C(S^2_{\mathbb{R}\mathcal{T}})}}$, for $i=0,1,2$, and by setting
\begin{gather*}
\alpha_0(1_{\mathcal{T}}\otimes u) :=\big(1_{{\mathcal{T}}}\otimes u,\phi_1\otimes 1_{C({\mathbb Z_2})},\phi_1\otimes 1_{C({\mathbb Z_2})}\big),
\\
\alpha_1(1_{\mathcal{T}}\otimes u) :=\big(\phi_1\otimes 1_{C({\mathbb Z_2})},1_{{\mathcal{T}}}\otimes u,\phi_2\otimes 1_{C({\mathbb Z_2})}\big),
\\
\alpha_2(1_{\mathcal{T}}\otimes u) :=\big(\phi_2\otimes 1_{C({\mathbb Z_2})},\phi_2\otimes 1_{C({\mathbb Z_2})},1_{{\mathcal{T}}}\otimes u\big).
\end{gather*}
Incidentally, the above formulas were not guessed but derived using the degenerated version of the construction from
Theorem~\ref{main-split} in which the relevant parts of $\alpha^i_j$'s were obtained utilizing $\phi_1$ and $\phi_2$.
By luck, the corrections in which the splittings $\beta^i_j$ could have been used turned out to be unnecessary~-- hence
it was also unnecessary to derive $\beta^i_j$'s using the methods from Section~\ref{partitions-s}.

Let us denote for brevity $\gamma_i:=\alpha_i(1_{\mathcal{T}}\otimes u)$, as well as omit subscripts indicating the
algebra the unit elements belong to.
Let us note that because $u^2=1$ we have
\begin{gather*}
1-\alpha_1^2=\big(\big(1-\phi_1^2\big)\otimes1,0,\big(1-\phi^2_2\big)\otimes 1\big),
\qquad
1-\alpha_1^2=\big(\big(1-\phi_2^2\big)\otimes1,\big(1-\phi^2_2\big)\otimes1,0\big).
\end{gather*}
Then the straightforward application of the formula from Theorem~\ref{Main} yields
\begin{gather*}
\ell(u) :=\alpha_0\otimes\alpha_0\big(1-\alpha_1^2\big)\big(1-\alpha_2^2\big) +\alpha_1\otimes\alpha_1\big(1-\alpha_2^2\big) +\alpha_2\otimes\alpha_2
\\
\phantom{\ell(u)\;}
 =(1\otimes u,\phi_1\otimes1,\phi_1\otimes 1) \otimes\big(\big(1-\phi_1^2\big)\big(1-\phi_2^2\big)\otimes u,0,0\big)
\\
\hphantom{\ell(u):=}{}
+(\phi_1\otimes1,1\otimes u,\phi_2\otimes 1) \otimes\big(\phi_1\big(1-\phi_2^2\big)\otimes1,\big(1-\phi_2^2\big)\otimes u,0\big)
\\
\hphantom{\ell(u):=}{}
+\big(\phi_2\otimes 1_{C({\mathbb Z_2})},\phi_2\otimes 1_{C({\mathbb Z_2})},1_{{\mathcal{T}}}\otimes u\big)
\otimes\big(\phi_2\otimes 1_{C({\mathbb Z_2})},\phi_2\otimes 1_{C({\mathbb Z_2})},1_{{\mathcal{T}}}\otimes u\big).
\end{gather*}
Note the similarity of this formula to the formula~\eqref{ellform3} obtained using the other method in the previous
subsection.
This similarity is understandable, because (not excluding the possibility of some general link, as yet unexplored by the
author, between the two methods used) the common feature of both particular computations is that by construction, both
strong connection formulas were expressed using the limited set of elements:
$\phi_1,\phi_2,1_{\mathcal{T}}\in{\mathcal{T}}$, and $1_{C({\mathbb Z_2})}, u\in C({\mathbb Z_2})$.

We leave to the reader analogous computations as those at the end of the previous subsection, which prove that both left
and right legs of the above strong connection are linearly independent (when taken separately).

\subsection*{Acknowledgements}

The author is grateful to Piotr M.~Hajac for helpful discussions.
The author would also like to thank the referees for helpful suggestions.
This work was partially supported by the NCN-grant 2012/06/M/ST1/00169.

\vspace{-2mm}

\pdfbookmark[1]{References}{ref}
\LastPageEnding

\end{document}